\documentclass[twoside,11pt]{amsart}

\usepackage{amsmath,latexsym,amssymb,mathptm, times}
\input amssym.def
\input amssym
\input xypic
\input xy
\xyoption{all}
\def\demo{\noindent{\bf Proof. }}
\def\sqr#1#2{{\vcenter{\hrule height.#2pt
        \hbox{\vrule width.#2pt height#1pt \kern#1pt
                \vrule width.#2pt}
        \hrule height.#2pt}}}
\def\square{\mathchoice\sqr64\sqr64\sqr{4}3\sqr{3}3}
\def\QED{\hfill$\square$}

\def\tratto{\mbox{\rule{2mm}{.2mm}$\;\!$}}

\def\M{{\overline M}}
\def\m{{\mathfrak m}}

\def\p{{\mathfrak p}}
\def\ms{\medskip}
\def\bs{\bigskip}

\newtheorem{Theorem}{Theorem}[section]
\newtheorem{Lemma}[Theorem]{Lemma}
\newtheorem{Corollary}[Theorem]{Corollary}
\newtheorem{Proposition}[Theorem]{Proposition}

\newtheorem{Notation and Discussion}[Theorem]{Notation and Discussion}
\newtheorem{Assumptions and Discussion}[Theorem]{Assumptions and Discussion}

\newtheorem{Remark}[Theorem]{Remark}
\newtheorem{Example}[Theorem]{Example}
\newtheorem{Definition}[Theorem]{Definition}

\begin{document}

\baselineskip=16pt

\title[$j$-multiplicity and depth of associated graded modules]
{\Large\bf $j$-multiplicity and depth of associated graded modules}

\author[C. Polini and Y. Xie]
{Claudia Polini \and Yu Xie}

\thanks{AMS 2010 {\em Mathematics Subject Classification}.
Primary 13A30; Secondary 13H15, 13B22, 13C14, 13C15, 13C40.}

\thanks{The first author was partially supported by the
NSF and the NSA}

\address{Department of Mathematics, University of Notre Dame,
Notre Dame, Indiana 46556} \email{cpolini@nd.edu}

\address{Department of Mathematics,University of Notre Dame,
Notre Dame, Indiana 46556} \email{yxie@nd.edu}

\vspace{-0.1in}

\begin{abstract}
Let $R$ be  a Noetherian local ring. We define the  minimal $j$-multiplicity and almost minimal $j$-multiplicity of an arbitrary $R$-ideal on any finite $R$-module.  For any ideal  $I$ with minimal $j$-multiplicity or almost minimal $j$-multiplicity on a Cohen-Macaulay module $M$, we
prove that under some residual assumptions, the associated graded module ${\rm gr}_I(M)$ is Cohen-Macaulay  or almost Cohen-Macaulay, respectively. Our work generalizes the
results for minimal multiplicity and almost minimal multiplicity
obtained by  Sally, Rossi, Valla, Wang, Huckaba, Elias, Corso, Polini, and VazPinto.
\end{abstract}

\maketitle

\vspace{-0.2in}

\section{Introduction}

In this paper we investigate the behavior of the depth of the
associated graded ring ${\rm gr}_I(R)$ of an ideal $I$ in a
Noetherian local ring $(R, \m)$ in terms of conditions on the
$j$-multiplicity of $I$. The associated graded ring of $I$ is an
algebraic construction whose projective scheme represents the
exceptional fiber of the blowup of a variety along a subvariety. Its
arithmetical properties, like its depth, provide useful information,
for instance, on the cohomology groups of the blowup. For an
$\m$-primary ideal $I$, the interplay between the Hilbert polynomial
of $I$, and more precisely its Hilbert coefficients, and the depth
of the associated graded ring has been widely investigated. This
line of study has its roots in the pioneering work of Sally.
The idea is that extremal values of the Hilbert coefficients, most
notably of the multiplicity of $I$,  yield high depth of the
associated graded ring and, conversely, good depth properties encode
information about all the Hilbert coefficients, such as their
positivity. The problem arises when one considers ideals which are
not $\m$-primary, because their Hilbert function is not defined,
thus there is no numerical information on Hilbert coefficients
available to study the Cohen-Macaulayness of ${\rm gr}_I(R)$. To
remedy the lack of this tool, in this paper we propose to use the
notion of $j$-multiplicity. The $j$-multiplicity was developed as a
generalization of the Hilbert-Samuel multiplicity to arbitrary
ideals. It was first introduced by Achilles and Manaresi in
1993 and, since then, it has been frequently used by both
algebraists and geometers as an invariant to deal with improper
intersections and varieties with non isolated singularities
\cite{AM}.

In this introduction we will only discuss the case of associated graded rings,
%, restrict ourself to the case of ideals,
although in the rest of the paper we will treat associated graded modules.
% arbitrary modules.

Let $I$ be an $\m$-primary ideal. The Hilbert-Samuel function of $I$
is the numerical function $H_I(n)$ that measures the growth of the
length $\lambda(R/I^n)$ of the powers of $I$ for all $n \geq 1$. For
$n$ sufficiently large, the Hilbert-Samuel function is a polynomial
function in $n$ of degree $d$, the dimension of $R$. This is the
Hilbert-Samuel polynomial of $I$, whose coefficients $e_i(I)$, dubbed the
Hilbert coefficients of $I$, are uniquely determined by $I$. It is
well known that the normalized leading coefficient $e_0(I)$, the
multiplicity of $I$, detects integral dependence of $\m$-primary
ideals. The integral closure of $I$, for instance, can be
characterized as the largest ideal containing $I$ with the same
multiplicity $e_0$, when the ring is equidimensional and universally
catenary. Flenner and Manaresi were the first to use the
$j$-multiplicity to generalize this fundamental theorem of Rees to
arbitrary ideals \cite{FM1}.

In 1967 Abhyankar proved that the multiplicity  $e(R)=e_0(\m)$ of a
$d$-dimensional Cohen-Macaulay local ring is bounded below by
$\mu(\m)-d+1$, where $\mu(\m)$ is the embedding dimension of $R$
\cite{AB}. Rings for which $e_0=\mu(\m)-d+1$ have since then been
called rings of minimal multiplicity. In the case of minimal
multiplicity, Sally had shown in \cite{S1} that the associated
graded ring of $\m$ is always Cohen-Macaulay. Even if the
multiplicity is almost minimal, the associated graded ring is
Gorenstein provided the ambient ring is Gorenstein \cite{S2}.
Unfortunately, for arbitrary Cohen-Macaulay rings of almost minimal
multiplicity the Cohen-Macaulayness of ${\rm gr}_\m(R)$ fails to
hold, the exceptions being Cohen-Macaulay local rings of maximal
type \cite{S3}. Based on this result, Sally conjectured that if the
multiplicity of $R$ is almost minimal then the depth of the
associated graded ring is almost maximal, i.e., it is at least
$d-1$. This conjecture was proved almost twenty years later by
Rossi and Valla \cite{RV1}, and independently by Wang
\cite{W}. In recent years there have been many generalizations of
these results to $\m$-primary ideals and modules with $I$-adic
filtrations, where $I$ is an ideal of definition, a condition that
is required to define the Hilbert-Samuel multiplicity (see for
example \cite{RV2}, \cite{H}, \cite{CPV}, \cite{E}, \cite{R1},  \cite{P1}, \cite{RV}).

The kind of generalization we accomplish in this article is much
harder. We investigate the depth of the associated graded ring of an
arbitrary ideal using the $j$-multiplicity introduced by Achilles
and Manaresi \cite{AM} and further studied in \cite{FOV}, \cite{FM1}, \cite{ FM2} and \cite{NU}.  Let $(R, \m)$ be a Noetherian local ring of dimension $d>0$
and $I$ an $R$-ideal. One can assign generalized Hilbert
coefficients $j_i(F)$ to every ideal filtration $F$ whose Rees
algebra is finite over $R[It]=\bigoplus_{i\geq 0} I^i$ in the following way:
let $A$ be the associated graded ring of $F$, and denote by
$\Gamma_{\m}(A) =H_{\m}^{0}(A)$ the submodule of elements supported
on $\m$. Since $\Gamma_{\m}(A)$ is annihilated by a large power of
$\m$, it is a finite graded module over ${\rm gr}_I(R) \otimes
R/{\m}^k$ for some $k$, hence its Hilbert polynomial is well
defined: $$P(n)= \sum_{i=0}^{d-1} (-1)^i j_i(F)
{{n+d-i-1}\choose{d-i-1}}.$$ We call $P(n)$ the {\it generalized
Hilbert polynomial} of $F$. The generalized Hilbert coefficients of
the filtration $\{I^n\}_{n \in {\mathbb N}}$ will simply be denoted
by $j_i(I)$. Notice that $j_0(I)$ coincides with the
$j$-multiplicity defined by Achilles and Manaresi in \cite{AM}.
Furthermore in the $\m$-primary case $j_i(I)=e_i(I)$, so our
definition recovers the standard one.

In Section 2, we prove a lower bound for the $j$-multiplicity of
any ideal $I$. The definition of ideals with minimal
$j$-multiplicity is thus immediate. In Section 3, under certain Artin-Nagata
condition, we prove that for any ideal with minimal
$j$-multiplicity, the associated graded ring
is Cohen-Macaulay. Furthermore if the ambient ring is Gorenstein, then the associated graded ring is Gorenstein as well, which generalizes completely Sally's results.
%The technical novelty is a powerful combination
%of the methods used in the $\m$-primary case with tools proper to
%residual intersection theory.
Finally, in Section 4, we deal with ideals  with almost minimal
$j$-multiplicity. We prove that, under the same residual assumptions, the
associated graded ring is almost
Cohen-Macaulay.
%Finally, in Section 4 we discuss the generalized
%Hilbert polynomials $P(n)$ for Cohen-Macaulay modules with minimal
%$j$-multiplicity or almost minimal $j$-multiplicity. Remarkably, the
%(almost) minimal $j$-multiplicity forces the shape of the
%generalized Hilbert polynomial and therefore determines the value of
%all the other generalized Hilbert coefficients.
The technical
novelty is a powerful combination of the methods used in the
$\m$-primary case with tools proper to residual intersection theory.
This result can be viewed as a positive answer to Sally's conjecture
for arbitrary ideals.

%By using this tool, one should also expect to extend the
%multiplicity theory of modules  with good ideal filtrations to
%arbitrary ideal filtrations.
%Despite of the fact that the Hilbert functions, the Hilbert
%coefficients as well as the depth of the associated graded modules
%are well understood for  Cohen-Macaulay modules with minimal
%Hilbert-Samuel multiplicity or almost minimal Hilbert-Samuel
%multiplicity,   very  little is known about what they could be for
%modules with the $j$-multiplicity.

%In this paper, we generalize the existing results of minimal multiplicity and
%almost minimal multiplicity for modules with good ideal filtrations to modules with arbitrary ideal filtrations.

\section{Minimal $j$-multiplicity}

\ms

%Let $(R,\m)$ be a Noetherian local ring and $I$ an $R$-ideal. Let $M$
%be a finite $R$-module    The ideal $I$ is called an
%{\it ideal of definition} for $M$ if $\lambda(M/IM)<\infty$, where
%$\lambda(M/IM)$ denotes the length of $M/IM$.  Since ${\rm dim}_G\, T={\rm dim}_R\, M=d$,  $T$ has a
%Hilbert function which is eventually a polynomial of degree $d-1$
%and the {\it Hilbert-Samuel multiplicity} of $M$ with respect to $I$
%is defined as
%$$%
%e(I, M)=(d-1)! {\it \mathop {\lim }\limits_{t \to \infty
%}}\frac{\lambda(I^tM/I^{t+1}M)}{t^{d-1}}.
%$$

%If $I$ is an $m$-primary ideal, $e(I, R)$ is the Hilbert-Samuel
%multiplicity of the ideal $I$ and $e(m, R)=e(R)$ the multiplicity
%of the local ring $R$.

%DEFINITION OF GOOD I-FILTRATION on M

In this section we first prove a lower bound for the
$j$-multiplicity; this bound leads to a notion of
minimal $j$-multiplicity.

We start by fixing the notation that will be used throughout the
paper. We first recall the definition of $j$-multiplicity according
to \cite{AM} and \cite{NU} .

Let $(R, \m)$ be a Noetherian local ring, $I$ an arbitrary
$R$-ideal, and $M$ a finite $R$-module of dimension $d$. The  {\it
$I$-adic filtration} of $M$ is a collection of submodules $\{I^jM\}_{j\geq 0}$. Let $G={\rm
gr}_I(R)=\oplus_{j=0}^{\infty}I^j/I^{j+1}$ be the {\it associated
graded ring} of $I$ and $T={\rm gr}_I (M)=\oplus_{j=0}^{\infty}
I^jM/I^{j+1}M$ the {\it associated graded module}  of the filtration $\{I^jM\}_{j\geq 0}$. Notice that $T$ is a finite graded module over the
graded ring $G$. In general the homogeneous components of  $T$ may
not have  finite length, thus we consider the $T$-submodule of
elements supported on $\m$,  $W=\Gamma_{\m}(T)=0:_T
{\m}^{\infty}=\oplus_{j=0}^{\infty} \Gamma_{\m}\,(I^jM/I^{j+1}M)$.
Since $W$ is annihilated by a large power of $\m$, it is a finite
graded module over ${\rm gr}_I(R) \otimes R/{\m}^k$ for some $k$,
hence its Hilbert polynomial $P(n)$ is well defined. Notice  that
 ${\rm dim}_G\,W\leq {\rm dim}_G\,T=d$, thus  $P(n)$
 has degree at most $d-1$. The $j${\it-multiplicity} of $I$ on
 $M$ is the normalized coefficient of $P(n)$ in degree $d-1$,
$$
j(I, M)=(d-1)! {\it \mathop {\lim }\limits_{t \to \infty
}}\frac{\lambda(\Gamma_\m(I^tM/I^{t+1}M))}{t^{d-1}}.
$$

Recall that the Krull dimension of the special fiber module
$T/{\m}T$ is called the
 {\it analytic spread} of $I$ on $M$ and is denoted by $\ell(I, M)$. In general,
 ${\rm dim}_G\,W \leq \ell(I,M) \leq d$
 and equalities hold if $\ell(I,M)=d$. Therefore $j(I, M)\neq 0$ if and only if $\ell(I, M)=d$ \cite[2.1]{NU}.

If $M/IM$ has finite length, the ideal $I$ is said to be  an {\it
ideal of definition} for $M$. In this case each homogeneous
component of $T$ has finite length, thus $W=T$ and the
$j$-multiplicity coincides with the usual multiplicity.

An element  $x\in I$ is said to be a  {\it superficial element}
for $I$ on $M$ if there exists a non-negative integer $c$ such that
$$
(I^{t+1}M:_M x)\cap I^cM=I^tM.
$$
A sequence of elements $x_1, \ldots, x_s$ in $I$ is a {\it superficial sequence}
for $I$ on $M$ if $x_i$ is superficial for $I$ on $M/(x_1,
\ldots, x_{i-1})M$ for $i=1, \ldots, s$.  This
notion, originally introduced by Zariski and Samuel, plays a
significant role in the study of Hilbert functions because it allows
to reduce the problems to lower dimensional ones. Notice that if $M$
has positive depth then every superficial element is regular on $M$.
If $d ={\rm dim} \, M \geq 1$ then a superficial element has always
order one, i.e., $x \in I/I^2$. Thus, in this case, the definition of
superficial elements coincides with the definition of homogeneous
filter regular elements used in the study of $j$-multiplicity (see
\cite{AM}, \cite{X}, and \cite{NU} for instance). More precisely, an element
$x$ is superficial for $I$ on $M$ if and only if $x^*$,
the image of $x$ in $I/I^2$, is {\it filter-regular} for $G_{+}$ on  $T$.

For an ideal $J\subseteq I$, one says that $J$ is a {\it reduction}
of $I$ on $M$ if $JI^tM=I^{t+1}M$ for some non negative
integer $t$. A {\it minimal} reduction is a reduction which is
minimal with respect to inclusion. Minimal reductions always exist
and, if $R$ has infinite residue field, the minimal number of
generators of any minimal reduction $J$ of $I$ on $M$
equals the analytic spread $\ell(I,M)$. Furthermore, a minimal
generating set of $J$ can be chosen to be a superficial sequence for
$I$ on $M$ \cite[3.1]{T1}. The least integer $t$ with
$JI^tM=I^{t+1}M$ is called the {\it reduction number} of $I$ on $M$
with respect to $J$ and denoted by  $r_J(I,M)$. One then defines the
{\it reduction number} $r(I,M)$ of $I$ on $M$ to be the least
$r_J(I,M)$, where $J$ varies over all minimal reductions of $I$ on $M$.

Let $I=(a_1, \ldots, a_n)$ and write $x_i=\sum_{j=1}^n
\lambda_{ij}a_j$ for $1\leq i\leq s$ and $(\lambda_{ij})\in R^{s
n}$.  The elements  $x_1, \ldots, x_s$ form a {\it sequence of
general elements} in $I$ (equivalently  $x_1, \ldots, x_s$ are {\it
general} in $I$) if there exists a dense open subset $U$ of $k^{sn}$
such that the image $(\overline{\lambda_{ij}})\in U$. In this case
we call the  $(\lambda_{ij})$ general elements in $ R^{s n}$. When
$s=1$ we say that $x=x_1$ is {\it general} in $I$. Observe that the
notion of  general elements of $I$ is more restrictive than the one
of {\it  sequentially general} elements. The latter means that for
every $i$ with $1 \leq i \leq s$ and every fixed $x_1, \ldots,
x_{i-1}$ the element $x_i$ is general in $I$.

The notion of general elements is a fundamental tool for our study
as they are always a superficial sequence for $I$ on
$M$ \cite[2.5]{X}; they generate a minimal reduction whose reduction number
$r_J(I, M)$ coincides with the reduction number $r(I,M)$ of $I$ on
$M$ if $s=\ell(I, M)$ (see \cite[2.2]{T2} and \cite[8.6.6]{SH}); and they form a {\it super-reduction} in the sense of \cite{AM} whenever
$s=\ell(I, M)=d={\rm dim}_R\,M$ (see \cite[2.5]{X}).  Furthermore, one can compute the $j$-multiplicity using general
elements and obtain a lower bound from it as the next proposition shows.

\begin{Proposition}\label{lowerbound}
Let $(R, \m)$ be a Noetherian local ring with infinite residue field
$k$. Let  $M$ be a finite $R$-module and $I$ an $R$-ideal with analytic spread $\ell(I, M)=d={\rm dim}_R\,M$. Then for general elements $x_1, \ldots, x_d$ in $I$, we have
\begin{itemize}

\item[(a)] the $j$-multiplicity of $I$ on $M$ is
\vskip -.3in
$$
j(I, M)=e(I, M/((x_1,\ldots,x_{d-1})M:_{M} I^{\infty}) )=
\lambda(M/((x_1,\ldots,x_{d-1})M:_{M} I^{\infty}+x_dM));
$$

\item[(b)]  $j(I,M) \geq  \lambda (I\overline{M}/I^2\overline{M})$ where $\overline{M}= M/((x_1,\ldots,x_{d-1})M:_M I^{\infty})$.
\end{itemize}
\end{Proposition}

\demo By \cite[2.5]{X}, there exist general  elements $x_1, \ldots, x_d$ in $I$ which form  a  super-reduction in the sense of \cite{AM}. Now proceed as in the proof of \cite[3.8]{AM} or \cite[3.6]{NU} to obtain the desired formula of part (a).
For part (b), as
$\overline{M}$ is  a one-dimensional Cohen-Macaulay
module and $I$ is an ideal of definition for $\overline{M}$, we have
$j(I, M)=e(I, \overline{M})=
\lambda(I\overline{M}/I^2\overline{M})+\lambda(I^2\overline{M}/x_d
I\overline{M})\geq \lambda (I\overline{M}/I^2\overline{M})$, where
the second equality follows from \cite[Corollary 2.1]{RV}.
\QED

\medskip

Once we show that $\lambda(I\overline{M}/I^2\overline{M})$ and
$\lambda(I^2\overline{M}/x_d I\overline{M})$ do not depend on the
choice of the general elements $x_1, \ldots, x_d$ in $I$, we will obtain
the desired lower bound for the $j$-multiplicity.

\begin{Lemma}\label{generic}
Let $(R, \m)$ be a Noetherian local ring with infinite residue field
$k$. Let $R^{\prime}=R[z_1, \ldots,z_t]$ be a polynomial ring over
$R$. Let $M^{\prime} \supseteq N^{\prime}$ be two finite
$R^{\prime}$-modules with
$\lambda_{R^{\prime}_{\m R^{\prime}}}(M^{\prime}_{\m R^{\prime}}/N^{\prime}_{\m R^{\prime}})=s$.
If $(\Lambda)=(\lambda_1, \ldots, \lambda_t)$ is a vector in
$R^{t}$, write $(\overline{\Lambda})$ for its image in $k^t$ and
$\pi(\cdot)$ for  the evaluation
map sending $z_i$ to $\lambda_i$. Then there exists a dense open
subset $U$ of $k^t$ such that
$\lambda_{R}(\pi(M^{\prime})/\pi(N^{\prime}))\leq s$ whenever
$(\overline{\Lambda})\in U$.
\end{Lemma}

\demo If $s=\infty$ then we are done. So we may assume $s<\infty$ and let $M^{\prime}=M_0^{\prime}  \supset  M_1^{\prime}\supset \ldots \supset M_s^{\prime}=N^{\prime}$
be a filtration such that $(M_{l-1}^{\prime}/M_{l}^{\prime})_{\m R^{\prime}}\cong k(z_1,\ldots,z_t)$
for $1\leq l\leq s$.
We use induction on $s$ to prove the lemma. When $s=0$, i.e.,
$\lambda_{R^{\prime}_{\m R^{\prime}}}(M^{\prime}_{\m R^{\prime}}/N^{\prime}_{\m R^{\prime}})=0$,
there exists a polynomial $f\in R^{\prime}\setminus \m R^{\prime}$ such that $
fM^{\prime}\subseteq N^{\prime}$. Let $\overline{f}$ be the image of $f$ in
$k[z_1, \ldots, z_t]$ and notice that $\overline{f}\neq 0$. Thus $U=D(\overline{f})$
is a dense open subset of $k^t$. If $(\overline{\Lambda})\in U$ then $f(\Lambda)$ is a unit
in $R$. Thus $f(\Lambda)\pi(M^{\prime})=\pi(fM^{\prime})\subseteq \pi(N^{\prime})$ implies
$\pi(M^{\prime})\subseteq \pi(N^{\prime})$.
Now assume the lemma holds for $s-1$, i.e., there exists a dense open subset
$U_1\subseteq k^t$ such that $\lambda_{R}(\pi(M^{\prime})/\pi(M_{s-1}^{\prime}))\leq s-1$
whenever $(\overline{\Lambda})\in U_1$.
Since $(M_{s-1}^{\prime}/M_{s}^{\prime})_{\m R^{\prime}}\cong k(z_1,\ldots,z_t)$,
there exists $b^{\prime}\in M_{s-1}^{\prime}$ and a polynomial $f\in R^{\prime}\setminus \m R^{\prime}$
so that  $fM_{s-1}^{\prime}\subseteq R^{\prime}b^{\prime}+M_{s}^{\prime}$ and
$f\m b^{\prime}\in M_{s}^{\prime}$.
%and $fb^{\prime}$ is not in $ M_{s}^{\prime}$.
Also notice that  the image $\overline{f}$ of $f$ in $k[z_1, \ldots, z_t]$ is not zero
and $U_2=D(\overline{f})$ is a dense open subset of $k^t$. Let $U=U_1\cap U_2$. Whenever $(\overline{\Lambda})\in U$,
 $f(\Lambda)$ is a unit in $R$. Thus $\pi(M_{s-1}^{\prime})\subseteq R\pi(b^{\prime})+\pi(M_{s}^{\prime})$
 and  $\m \pi(b^{\prime})\in \pi(M_{s}^{\prime})$.
 %and $\pi(b^{\prime})$ is not in $ \pi(M_{s}^{\prime})$.
 Therefore $\lambda_{R}(\pi(M_{s-1}^{\prime})/\pi(M_s^{\prime}))\leq 1$ and thus we obtain $\lambda_{R}(\pi(M^{\prime})/\pi(N^{\prime}))=\lambda_{R}(\pi(M^{\prime})/\pi(M_{s-1}^{\prime}))+ \lambda_{R}(\pi(M_{s-1}^{\prime})/\pi(M_s^{\prime}))\leq s$ for every  $(\overline{\Lambda})\in U$. \QED

\ms

\begin{Lemma}\label{general}
Let $(R,\m)$ be a Noetherian local ring with infinite residue field $k$.  Let  $M$ be a finite $R$-module and  $I$ an $R$-ideal with analytic spread $\ell(I, M)=d={\rm dim}_R\,M$. For $x_1, \ldots, x_d$ general elements in $I$, write $\overline{M}=M/((x_1,\ldots,x_{d-1})M:_M I^{\infty})$, then the lengths  $\lambda(I\overline{M}/I^2\overline{M})$ and $\lambda(I^2\overline{M}/x_d I\overline{M})$ are independent of  $x_1, \ldots, x_d$.
\end{Lemma}

\demo
Let $a_1, \ldots, a_n$ be a set of generators of $I$ and $Z=(z_{ij})$ be $d\times n$ variables. Write  ${R^{\prime}}=R[Z]$, $x_i^{\prime}=\sum_{j=1}^{n}z_{ij}a_j$ for $1\leq i\leq d$, and $M^{\prime}=M\otimes_R {R^{\prime}}$.  Let $\overline{M^{\prime}}=M^{\prime}/((x_1^{\prime},\ldots,x_{d-1}^{\prime})M^{\prime}:_{M^{\prime}} {IR^{\prime}}^{\infty})$, by Proposition~\ref{lowerbound}(a) and the proof of  Proposition~\ref{lowerbound}(b) (see also \cite[3.8]{AM}, \cite[3.6]{NU}  and \cite[Corollary 2.1]{RV})  we have
$$j(I, M)=j(I{R^{\prime}_{\m R^{\prime}}}, M^{\prime}_{\m R^{\prime}})
%=\lambda(M^{\prime}_{\m R^{\prime}}/[(x_1^{\prime},\ldots,x_{d-1}^{\prime})
%M^{\prime}_{\m R^{\prime}}:_{M^{\prime}_{\m R^{\prime}}} {(IR^{\prime}_{\m R^{\prime}}})^{\infty}+x_d^{\prime}M^{\prime}_{\m R^{\prime}}])
=\lambda(I\overline{M^{\prime}}_{\m R^{\prime}}/I^2\overline{M^{\prime}}_{\m R^{\prime}})+
\lambda(I^2\overline{M^{\prime}}_{\m R^{\prime}}/x_d^{\prime}I\overline{M^{\prime}}_{\m R^{\prime}}).
$$

For general elements $(\Lambda)=(\lambda_{ij})\in R^{dn}$, write $\pi(\cdot)$ for the  evaluation map sending $z_{ij}$ to $\lambda_{ij}$. Observe $\pi(IM^{\prime})= IM$, $\pi(I^2M^{\prime})=I^2M$, $\pi(x_d^{\prime}IM^{\prime})= x_d IM$, and clearly
$$\pi((x_1^{\prime},\ldots,x_{d-1}^{\prime})M^{\prime}:_{M^{\prime}} I{R^{\prime}}^{\infty})\subseteq (x_1,\ldots,x_{d-1})M:_{M} I^{\infty}.$$
Putting this together with Lemma \ref{generic}, we obtain
\begin{eqnarray*}
 \lambda(I\overline{M^{\prime}}_{\m R^{\prime}}/I^2\overline{M^{\prime}}_{\m R^{\prime}}) &=& \lambda(IM^{\prime}_{\m R^{\prime}}/[(x_1^{\prime},\ldots,x_{d-1}^{\prime})M^{\prime}
_{\m R^{\prime}}:_{IM^{\prime}_{\m R^{\prime}}} (I{R^{\prime}}_{\m R^{\prime}})^{\infty}+I^2M^{\prime}_{\m R^{\prime}}])       \\
   &\ge&\lambda(IM/[\pi((x_1^{\prime},\ldots,x_{d-1}^{\prime})M^{\prime}:_{IM^{\prime}} I{R^{\prime}}^{\infty})+I^2M])  \\
   &\ge& \lambda(IM/[(x_1,\ldots,x_{d-1})M:_{IM} I^{\infty}+I^2M])=\lambda(I\overline{M}/I^2\overline{M}).
\end{eqnarray*}

In the same way we have,
\begin{eqnarray*}\lambda(I^2\overline{M^{\prime}}_{\m R^{\prime}}/x_d^{\prime}I\overline{M^{\prime}}_{\m R^{\prime}})
&=&\lambda(I^2M^{\prime}_{\m R^{\prime}}/[(x_1^{\prime},\ldots,x_{d-1}^{\prime})M^{\prime}_{\m R^{\prime}}
:_{I^2M^{\prime}_{\m R^{\prime}}} (I{R^{\prime}}_{\m R^{\prime}})^{\infty}+x_d^{\prime}IM^{\prime}_{\m R^{\prime}}])\\
&\geq& \lambda(I^2M/[(x_1,\ldots,x_{d-1})M:_{I^2M} I^{\infty}+x_dIM])=
\lambda(I^2\overline{M}/x_dI\overline{M}).
\end{eqnarray*}

By Proposition~\ref{lowerbound}(a) and the proof of Proposition~\ref{lowerbound}(b), the $j$-multiplicity is given by the sum of $\lambda(I\overline{M}/I^2\overline{M})$ and $\lambda(I^2\overline{M}/x_dI\overline{M})$, thus
\begin{eqnarray*}
j(I, M)&=& \lambda(I\overline{M}/I^2\overline{M})+\lambda(I^2\overline{M}/x_dI\overline{M})\\
&\le&\lambda(I\overline{M^{\prime}}_{\m R^{\prime}}/I^2\overline{M^{\prime}}_{\m R^{\prime}})+
\lambda(I^2\overline{M^{\prime}}_{\m R^{\prime}}/x_d^{\prime}I\overline{M^{\prime}}_{\m R^{\prime}})\\
&=& j(I,M).
\end{eqnarray*}
In turn this forces the equalities
\begin{eqnarray*}
\lambda(I\overline{M}/I^2\overline{M})&=&
\lambda(I\overline{M^{\prime}}_{\m R^{\prime}}/I^2\overline{M^{\prime}}_{\m R^{\prime}}), \\
\lambda(I^2\overline{M}/x_dI\overline{M})&=&
\lambda(I^2\overline{M^{\prime}}_{\m R^{\prime}}/x_d^{\prime}I\overline{M^{\prime}}_{\m R^{\prime}}),
\end{eqnarray*} and therefore shows the independence of these lengths from the general elements $x_1, \ldots, x_d$.
\QED

\ms

Because of Proposition~\ref{lowerbound} and Lemma~\ref{general}, we can now give the definition of {\it minimal $j$-multiplicity of $I$ on $M$} which is the analogue of minimal multiplicity  \cite{RV}.

\begin{Definition}\label{j-multiplicity}
{\rm   Let $M$ be a finite module of dimension $d$ over a Noetherian local ring $R$ and
  $I$ an $R$-ideal with analytic spread $\ell(I, M)=d$. We say that $I$ has {\it minimal $j$-multiplicity on  $M$} if  $j(I, M)= \lambda(I\overline{M}/I^2\overline{M})$, where $\overline{M}= M/((x_1,\ldots,x_{d-1})M:_M I^{\infty})$ and $x_1, \ldots, x_{d-1}$ are general in $I$.}
\end{Definition}

Notice that   $I$ has  minimal $j$-multiplicity on  $M$ if and only if  $\lambda(I^2\overline{M}/x_dI\overline{M})$ is zero, or equivalently, if and only if $x_d$ generates a reduction of $I$ on $\M$ with reduction number one. The next lemma shows that for an ideal the assumption of having minimal $j$-multiplicity on $M$ is quite strict.
Indeed, if $I$ has  minimal $j$-multiplicity on  $M$ then the Hilbert function of $I$ on $\M$ is rigid, i.e., the value of the multiplicity determines the  Hilbert function. We remark that results of this kind are really surprising since the multiplicity is just one of the Hilbert coefficients and, in turn, the Hilbert coefficients give only partial information on the Hilbert polynomial which gives only asymptotic information on the Hilbert function.

\ms
\begin{Corollary}\label{rigidity} Let $R$, $I$, $M$ and $\overline{M}$ be as in Definition \ref{j-multiplicity}, then the $j$-multiplicity of $I$ on $M$ and on $I^tM$ coincides for every $t \ge 0$, i.e., $j(I, M)=j(I, I^tM)$ for every $t\geq 0$. Furthermore, if $I$ has  minimal $j$-multiplicity on  $M$  then $j(I, M)=\lambda(I^t\overline{M}/I^{t+1}\overline{M})$ for every $t\geq 1$.
\end{Corollary}

\demo
Observe that $\ell(I, M)={\rm dim} \, T/\m T={\rm dim} \, T_{\ge t}/\m T_{\ge t}=\ell(I, I^tM)$ and if $x_1, \ldots, x_d$ are general elements of $I$ on $M$, they are also general on $I^tM$ for all $t \ge 0$. The first assertion follows because $(x_1,\ldots,x_{d-1})M:_{IM} I^{\infty}= (x_1,\ldots,x_{d-1})IM:_{IM} I^{\infty},\ $ thus \begin{eqnarray*}
j(I, M)&=&\lambda(IM/((x_1,\ldots,x_{d-1})M:_{IM} I^{\infty}+x_dIM))\\
&=&\lambda(IM/((x_1,\ldots,x_{d-1})IM:_{IM}I^{\infty}+x_dIM) )\\
&=&j(I, IM).
\end{eqnarray*}
Now continue the above process, we will obtain $j(I, M)=j(I, I^tM)$ for every $t\geq 0$.

For the second assertion assume  that $I$ has  minimal $j$-multiplicity on  $M$, or equivalently, $I^2\M=x_dI\M$. This gives  $I^{t+2}\overline{M}=x_d I^{t+1}\overline{M}$ for every $t\geq 0$. In turn, this implies that $x_d$ generates a reduction of $I$ on $I^t\M$ with reduction number one, thus $I$ has  minimal $j$-multiplicity on  $I^tM$ and therefore $j(I, I^tM)=\lambda(I\overline{I^t M}/I^{2}\overline{I^tM})=\lambda(I^{t+1}\overline{M}/I^{t+2}\overline{M}) .$ Hence $j(I,M)= j(I, I^tM)=\lambda(I^{t+1}\overline{M}/I^{t+2}\overline{M})$ for every $t \ge 0$.
\QED

\ms

\section{Cohen-Macaulayness of the associated graded module }

\ms

In this section we show that the associated graded
module of any filtration with minimal $j$-multiplicity is
Cohen-Macaulay, if the ideal has some residual properties. We start by describing the residual assumptions that are needed to prove the main theorem.

Let $M$ be a   finite faithful module of dimension $d$  over a
Noetherian local ring $R$. Let $I$ be an $R$-ideal. The ideal $I$ is
said to satisfy the {\it condition $G_s$} on  $M$ if for every
$\p\in {\rm Supp}_R(M/IM)$ with ${\rm dim}\,R_\p=t<  s$ the ideal
$I$ is generated by $t$ element on $M_\p$, i.e., $IM_\p=(x_1,
\ldots, x_{t})M_\p$ for some $x_1, \ldots, x_{t}$ in $I$. Write
$H=(x_1, \ldots, x_t)M:_M I$. If $IM_\p=(x_1, \ldots, x_{t})M_\p$
for every $\p\in {\rm Spec}(R)$ with ${\rm dim}\,R_\p\leq t-1$, then
$H$ is said to be a {\it  $t$-residual intersection } of $I$ on $M$.
Now let $H$ be a $t$-residual intersection of $I$ on $M$. If in
addition $IM_\p=(x_1, \ldots, x_{t})M_\p$ for every  $\p\in {\rm
Supp}_R(M/IM)$ with ${\rm dim}\,R_\p\leq t$, then $H$ is said to be
a {\it geometric $t$-residual intersection } of $I$ on $M$. If $M$
is not faithful, then we say that $I$  satisfies the {\it condition
$G_s$} on $M$ if $I\overline{R}$  satisfies the condition $G_s$ on
$M$, where $\overline{R}=R/{\rm Ann}\,M$. We say $H$ is a {\it
$t$-residual intersection } or {\it geometric $t$-residual
intersection } of $I$ on $M$ if $H$ is a $t$-residual intersection
or geometric $t$-residual intersection of $I\overline{R}$ on $M$
respectively.

\ms

The next two lemmas contain basic facts about residual intersections. The ideas are already presented in \cite[1.6 and 1.7]{U}.
The first lemma says that  the condition $G_{s}$  gives rise to residual intersections.

\begin{Lemma}\label{generalresidual}
Let $R$ be a catenary and equidimensional Noetherian local ring with infinite residue field. Let  $M$  be a finite $R$-module and $I$ an $R$-ideal satisfying condition $G_s$ on $M$. For general elements $x_1,  \ldots, x_s$ of $I$ on $M$, write $H_i=(x_1, \ldots, x_i)M:_M I $ for $0\leq i\leq s$, then:

$($a$)$ $H_i$ is an $i$-residual intersection of $I$ on $M$ for $0\leq i\leq s$.

$($b$)$ $H_i$ is a geometric $i$-residual intersection of $I$ on  $M$ for $0\leq i\leq s-1$.

\end{Lemma}

\demo  Let $a_1, \ldots, a_n$ be a set of generators of $I$ on $M$, i.e, $IM=(a_1, \ldots, a_n)M$, and $Z=(z_{ij})$ be  $s\times n$  variables. Write $R^{\prime}=R[Z]_{\m R[Z]}$ and  $M^{\prime}=M\otimes_R R^{\prime}$, where $\m$ is the maximal ideal of $R$. For $1\leq i\leq s$, let $x_i^{\prime}=\sum_{j=1}^nz_{ij}a_j$.
We first claim that

($\star$) if $\p\in {\rm Supp}(M/IM)$ such that   $IM_\p=(y_1, \ldots, y_l)M_\p$ for some $y_1, \ldots, y_l\, \in I$ where $1\leq l\leq s$, then  $IM_\p^{\prime}=(x_1^{\prime},\ldots, x_l^{\prime})M^{\prime}_\p$.

To prove it,  let $X_1, \ldots, X_n$ be variables over $R^{\prime}_\p/\p R^{\prime}_\p$. Set  $A^{\prime}=M^{\prime}_\p/\p M^{\prime}_\p [X_1, \ldots, X_n]$, $T^{\prime}_\p={\rm gr}_{IR^{\prime}_\p}(M^{\prime}_\p)$, and $F^{\prime}=T^{\prime}_\p/\p T^{\prime}_\p$. Define the map $\varphi: A^{\prime}\rightarrow F^{\prime}$ by sending $X_i$ to $a_i+\p I R^{\prime}_\p$. For $1\leq i \leq l$, write $y_i=\sum_{j=1}^n\lambda_{ij}a_j$, where $(\lambda_{ij})\in R^{ln}$. Denote with $^{\tratto}$ the images of elements of $R$ in $R^{\prime}_\p/\p R^{\prime}_\p$. Observe that the preimages of the $y_i's$ generate a vector space in  $[A^{\prime}]_1$ of dimension $n$, i.e., ${\rm dim}\ {\rm Span} \{\{b_i=\sum_{j=1}^n\overline{\lambda_{ij}}X_j\}_{ 1\leq i \leq l}, \  [{\rm ker}(\varphi)]_1\}=n$. To show that $IM_\p^{\prime}=(x_1^{\prime},\ldots, x_l^{\prime})M^{\prime}_\p$ it will be enough to show  that the set  $ \{\{b_i^{\prime}=\sum_{j=1}^nz_{ij}X_j \}_{ 1\leq i \leq l}, \  [{\rm ker}(\varphi)]_1\}$ spans also a $n$-dimensional vector space in $[A^{\prime}]_1$, i.e., ${\rm dim}\ {\rm Span}  \{\{b_i^{\prime}=\sum_{j=1}^nz_{ij}X_j\}_{ 1\leq i \leq l} , \  [{\rm ker}(\varphi)]_1\}=n$. Indeed,  define the map $$\phi: A=M_\p/\p M_\p [X_1, \ldots, X_n]\rightarrow F={{\rm gr}_{IR_{\p}}}(M_\p)/(\p{{\rm gr}_{IR_{\p}}}(M_\p))$$ by sending $X_i$ to $a_i+\p I R_\p$. Since the extension from $R$ to $R^{\prime}$ is flat, ${\rm ker}(\varphi)={\rm ker}(\phi)\otimes_R R^{\prime}$. Thus we can choose basis $b_i=\sum_{j=1}^n\overline{\lambda_{ij}}X_j$ in $[{\rm ker}(\phi)]_1$  where $\lambda_{ij}\in R$ and $l+1\leq i\leq t $ such that $${\rm dim}\ {\rm Span} (b_1, \ldots b_{l+t}) =n.$$ This forces  ${\rm dim}\ {\rm Span} (b_1^{\prime}, \ldots, b_l^{\prime}, b_{l+1}, \ldots, b_{l+t}) =n$ as well. If not, set $$\alpha={\rm dim}\ {\rm Span} (b_1^{\prime}, \ldots, b_l^{\prime}, b_{l+1}, \ldots, b_{l+t}) < n.$$ Observe that each $b_i^{\prime}$ is a linear combination of $X_j$ with coefficients polynomials in the variable $z_{ij}$ of degree at most one, and let $X$ be the $(l+t)\times n$ matrix obtained from these coefficients. Because $\alpha < n$ then all the $n\times n$ minors of $X$ vanish. When we specialize $z_{ij}$ to $\lambda_{ij}$ for $1\leq i\leq l$, then  all the $n\times n$ minors of the $(l+t)\times n$ matrix obtained from the coefficients of $b_1, \ldots, b_l, b_{l+1}, \ldots, b_{l+t}$ as linear combinations of $X_j$ vanish as well, which contradicts the fact that ${\rm dim}\ {\rm Span} (b_1, \ldots b_{l+t}) =n$.

\smallskip

Write $H_i^{\prime}=(x_1^{\prime}, \ldots, x_i^{\prime})M^{\prime}:_{M^{\prime}} IR^{\prime} $ for $0\leq i\leq s$. We claim that:

\noindent
($a^{\prime}$) $H_i^{\prime}$ is an  $i$-residual intersection of $IR^{\prime}$ on $M^{\prime}$ for $0\leq i\leq s$.

\noindent
 ($b^{\prime}$)  $H_i^{\prime}$ is a geometric $i$-residual intersection of $IR^{\prime}$ on $M^{\prime}$ for $0\leq i\leq s-1$.

\smallskip

 Now by factoring out ${\rm Ann}\,M$, we may assume $M$ is faithful. For part ($a^{\prime}$),  let $1\leq i\leq s$ and $\p^{\prime}\in {\rm Spec}(R^{\prime})$ with ${\rm ht}\,\p^{\prime}\leq i-1$. Then $\p^{\prime}=\p R^{\prime}$ for some $\p\in {\rm Spec}(R)$ with ${\rm ht}\,\p\leq i-1$. If $\p \not\in {\rm Supp}(M/IM)$, then $IM_{\p^\prime}^{\prime}=M_{\p^\prime}^{\prime}=(x_1^{\prime},\ldots, x_i^{\prime})M^{\prime}_{\p^\prime}$, where the last equality holds because $x_i^{\prime}R_{\p^{\prime}}^{\prime}=R_{\p^{\prime}}^{\prime}$. Therefore we may assume $\p \in {\rm Supp}(M/IM)$. Since $I$  satisfies condition $G_s$ on $M$ and $\p\in {\rm Spec}(R)$ with ${\rm ht}\,\p\leq i-1<s$, we have
 $IM_\p=(y_1, \ldots, y_{i-1})M_\p$ for some $y_1, \ldots, y_{i-1}\, \in I$. By ($\star$), this implies $IM_\p^{\prime}=(x_1^{\prime},\ldots, x_{i-1}^{\prime})M^{\prime}_\p$.  Part  ($b^{\prime}$) follows by employing the same argument.

\smallskip
Finally since for $0\leq i\leq s-1$,
 $H_i^{\prime}$ is a geometric  $i$-residual intersection of
 $IR^{\prime}$ on $M^{\prime}$ and $H_s^{\prime}$ is a   $s$-residual intersection of
 $IR^{\prime}$ on $M^{\prime}$, we have ${\rm ht}\,( (x_1^{\prime},\ldots, x_i^{\prime})M^{\prime}:_{R^{\prime}} IM^{\prime})\geq i$ for $0\leq i\leq s$ and ${\rm ht}\,((x_1^{\prime},\ldots, x_i^{\prime})M^{\prime}:_{R^{\prime}} IM^{\prime}+IM^{\prime}:_{R^{\prime}} M^{\prime})\geq
 i+1$ for $0\leq i\leq s-1$. Let $k=R/\m$, $(\Lambda)=(\lambda_{ij})\in
 R^{sn}$ and $(\overline{\Lambda})$  the image of  $(\Lambda)$ in $k^{sn}$.
 Write
 $\pi(\cdot)$ for  the evaluation map sending $z_{ij}$ to $\lambda_{ij}$.
 By \cite[3.1]{HHU} for
 all $i$, there exists a dense open subset $U$ of $k^{sn}$
 such that ${\rm ht}(\pi((x_1^{\prime},\ldots, x_i^{\prime})M^{\prime}:_{R^{\prime}} IM^{\prime}))
 \geq i$ and ${\rm ht}(\pi((x_1^{\prime},\ldots, x_i^{\prime})M^{\prime}:_{R^{\prime}} IM^{\prime}+IM^{\prime}:_{R^{\prime}} M^{\prime}))\geq i+1$ whenever $(\overline{\Lambda})\in U$.
Since $\pi((x_1^{\prime},\ldots, x_i^{\prime})M^{\prime}:_{R^{\prime}} IM^{\prime}))\subseteq (x_1, \ldots, x_i)M:_R IM$ and $\pi(IM^{\prime}:_{R^{\prime}} M^{\prime})\subseteq IM:_{R} M$, then $H_i$ is also a geometric  $i$-residual intersection of
 $I$ on $M$ for  $0\leq i\leq s-1$ and $H_s$ is a   $s$-residual intersection of
 $I$ on $M$.
\QED

\ms
Assume that $M$ is   Cohen-Macaulay. The ideal $I$ is said to have the {\it Artin-Nagata property} ${AN^-_t}$ on $M$ if for every $0\leq i\leq t$ and every geometric $i$-residual intersection $H$ of $I$ on $M, \,$  the module $M/H$ is  Cohen-Macaulay.
%If $M$ is not faithful, then we say that $I$ has the {\it Artin-Nagata property} ${AN^-_t}$ on $M$ if $I\overline{R}$  has the {\it Artin-Nagata property} ${AN^-_t}$ on $M$ where $\overline{R}=R/{\rm ann}\,M$.
In the next lemma we exhibit some basic facts about Artin-Nagata properties that will be useful in the proof of Theorem~\ref{minimal}.

\ms

\begin{Lemma}\label{moduleresidual}%$($see \cite[2.3, 2.4]{JU} or \cite[1.7]{U}$)$
Let $M$ be a Cohen-Macaulay module of dimension $d$ over a Noetherian local ring $R$ with infinite residue field. Let $I$ be an $R$-ideal with $\ell(I,M)=s$ satisfying $G_{s}$ and  $AN^-_{s-1}$ on $M$.  For general elements $x_1,  \ldots, x_s$ of $I$ on $M$, write $H_i=(x_1, \ldots, x_i)M:_M I $ for $0\leq i\leq s$, then:

$($a$)$ $x_{i+1}$ is regular on $M/H_i$ and $H_i=(x_1, \ldots, x_i)M:_M x_{i+1}$ \  if \   $0\leq i\leq s-1$.

$($b$)$ $M/H_i$ is unmixed of dimension $d-i$.

$($c$)$  ${\rm depth}\, M/(x_1, \ldots, x_i)M=d-i$.

$($d$)$  ${\rm depth}\,(M/(x_1, \ldots, x_i)IM)\geq {\rm min}\{d-i, {\rm depth}\,(M/IM)\}$.

$($e$)$   $(x_1, \ldots, x_i)M=H_i\cap IM $ \  if \   $0\leq i\leq s-1$.

%and ${\rm ht}\, \overline{I}=1$ if $g-1\leq i\leq s-1$, where $g={\rm ht}\,I$.

%$($e$)$ Let $p\in V(I)$, then $I_{\p}$ satisfies $AN^-_{s-1}$ with respect to $M_{\p}$.

$($f$)$ If ${\rm depth}\,(M/IM)\geq d-s+1$ then $$(x_1,\ldots,x_{s-1})M:_{I^2M}I^{\infty}=  (x_1,\ldots,x_{s-1})M:_{I^2M}I= (x_1,\ldots,x_{s-1})M:_{I^2M}x_s=(x_1,\ldots,x_{s-1})IM.$$

$($g$)$ Let $\overline{M}=M/H_0$ then $I $ satisfies $G_s$ and $AN^-_{s-1}$  on $\, \overline{M}$.
%(c) $(x_1, \ldots, x_i)M=IM\cap H_i$ if $0\leq i\leq s-1$, or if $i=s$ and $H$ is a geometric $s$-residual intersection.

%(d) ${\rm ht}I^{\prime}=1$ if $g-1\leq i\leq s-1$, or if $i=s$ and $H$ is a geometric $s$-residual intersection.

%(e) $R^{\prime}/I^{\prime}\cong R/I+H_i$ is Cohen-Macaulay if $0\leq i\leq s-1$ and $R/I$ is cohen-Macaulay.

%(f) $H_{i+1}^{\prime}=(a_{i+1}^{\prime}):I^{\prime} $ and $a_{i+1}^{\prime})$ is regular on $R^{\prime}$, if $0\leq i\leq s-1$.

\end{Lemma}

\demo By Lemma \ref{generalresidual}, for general elements $x_1,  \ldots, x_s$ of $I$ on $M$, the module $H_i$ is a geometric $i$-residual intersection of $I$ on $M$ for $0\leq i\leq s-1$ and $H_s$ is a $s$-residual intersection of $I$ on $M$. Thus parts  (a), (b), (c), (e) and (g) follow as in the proofs of  \cite[1.7]{U} and \cite[2.3, 2.4]{JU}.   To prove (d),
we use induction on $i$.
First when $i=0$,  clearly ${\rm depth}\,M\geq {\rm min}\{d, {\rm depth}\,(M/IM)\}$. Assume  ${\rm depth}\,(M/(x_1, \ldots, x_i)IM)\geq {\rm min}\{d-i, {\rm depth}\,(M/IM)\}$ for some $0\leq i< s$. By  (a) and (e), $((x_1, \ldots, x_i)M:_Mx_{i+1})\cap IM=(x_1, \ldots, x_i)M$. For \,$i+1$, we have
\begin{eqnarray*}
(x_1, \ldots, x_i)IM\cap x_{i+1}IM&=&x_{i+1}[((x_1, \ldots, x_i)IM:_Mx_{i+1})\cap IM] \\
&\subseteq& x_{i+1}[((x_1, \ldots, x_i)M:_M x_{i+1})\cap IM] \\
&=&x_{i+1}(x_1, \ldots, x_i)M\subseteq (x_1, \ldots, x_i)IM\cap x_{i+1}IM.
\end{eqnarray*}
Thus we obtain an exact sequence:
$$
0\rightarrow x_{i+1}(x_1, \ldots, x_i)M\rightarrow (x_1, \ldots, x_i)IM\oplus x_{i+1} IM \rightarrow (x_1, \ldots, x_{i+1})IM \rightarrow 0.
$$
The element $x_{i+1}$ is regular on $IM$ and therefore on $(x_1, \ldots, x_i)M$ because $(0_M:_M x_{i+1})\cap
(x_1, \ldots, x_i)M $
$\subseteq (0_M:_Mx_{i+1})\cap IM=0_M$, thus $x_{i+1}
(x_1, \ldots, x_i)M\cong (x_1, \ldots, x_i)M$ and $x_{i+1} IM\cong IM$ .
%It is clear when $i+1=d$. So we may assume $i+1<d$.
Therefore ${\rm depth}\,(x_{i+1}(x_1, \ldots, x_i)M)\geq {\rm min}\{d, d-i+1\}$ and ${\rm depth}\,(x_{i+1} IM)\geq {\rm min}\{d, {\rm depth}\,(M/IM)+1\}$. By induction hypothesis ${\rm depth}\,((x_1, \ldots, x_i)IM)\geq {\rm min}\{d-i+1, {\rm depth}\,(M/IM)+1\}$, thus the above exact sequence yields ${\rm depth}\,((x_1, \ldots, x_{i+1})IM)\geq {\rm min}\{ d-i, {\rm depth}\,(M/IM)+1 \}.$ We finally conclude that
${\rm depth}\,(M/(x_1, \ldots, x_{i+1})IM)\geq {\rm min}\{d-i-1, {\rm depth}\,(M/IM)\}.$

To show assertion (f), since $$(x_1,\ldots,x_{s-1})IM \subseteq (x_1,\ldots,x_{s-1})M:_{I^2M}I\subseteq (x_1,\ldots,x_{s-1})M:_{I^2M}I^{\infty},$$ it is enough to check the equality locally at every prime ideal ${\p}\in {\rm Ass}(M/(x_1,\ldots,x_{s-1})IM)$.
By   (d), ${\rm depth}\,(M/(x_1,\ldots,x_{s-1})IM)\geq d-s+1.$ Thus for every  ${\p}\in {\rm Ass}(M/(x_1,\ldots,x_{s-1})IM)$, ${\rm ht}\,{\p}\leq s-1$ and hence either $IM_{\p}=M_{\p}$ or $IM_{\p}=(x_1, \ldots, x_{s-1})M_{\p}$ since $H_{s-1}$ is a geometric $s-1$-residual intersection of $I$ on $M$. Therefore, if  $IM_{\p}=(x_1, \ldots, x_{s-1})M_{\p}$ then $$[(x_1,\ldots,x_{s-1})M:_{I^2M}I^{\infty}]_{\p}=[(x_1,\ldots,x_{s-1})M:_{I^2M}I]_{\p}=(x_1, \ldots, x_{s-1})IM_{\p}= I^2M_{\p}.$$ Otherwise by (e), we immediately obtain $(x_1, \ldots, x_{s-1})M_{\p}= (H_{s-1})_{\p}$, and the latter coincides with the module $[(x_1,\ldots,x_{s-1})M:_{I^2M}I]_{\p}=[(x_1,\ldots,x_{s-1})M:_{I^2M}I^{\infty}]_{\p}.$

The full assertion now follows from part (a).

\QED

The following lemma shows that  the presence of $G_{d}$ along with the Artin-Nagata condition $AN^-_{d-2}$ is actually sufficient to obtain $AN^-_d$.

\begin{Lemma}\label{Moduleresidual}$(${\rm see} \cite[1.9]{U}$)$
Let $M$ be a  Cohen-Macaulay module of dimension $d$ over a Noetherian local ring $R$ and  $I$  an $R$-ideal  satisfying $G_{d}$ and  $AN^-_{d-2}$ on $M$. Then
$I$ satisfies  $AN^-_{d}$ on $M$.
\end{Lemma}
\demo
The claim follows from Lemma \ref{moduleresidual} (b).
\QED

%CHECK THE REMARK
%The next remark affirms that  super reductions give rise to residual intersections. It is essentially taken by \cite[1.11]{U} and \cite[2.7]{JU}. We state here for convenience the module version. The proof is the same as in the case of ideals.

%\begin{Remark}\label{RI}
%Let $M$ be a Cohen-Macaulay module of dimension $d$ over a Noetherian local ring $R$ with infinite residue field.  Let $I=(a_1, \ldots,a_n)$ be an $R$-ideal with $\ell(I, M)=d$. Assume  ${\rm depth}\,(M/IM)\geq {\rm min}\{{\rm dim} (M/IM), 1\}$ and  $I$ satisfies $G_{d}$ and  $AN^-_{d-3}$ with respect to $M$.  If   $x_1, \ldots, x_d$ is a super-reduction of $I$
%with respect to $M$ then
% ${\rm ht}\, ({\rm ann}_R\,M/(x_1, \ldots, x_d)M:_M I)\geq d$.
 %\end{Remark}

 \ms

Now we show that, as in the $\m$-primary case, minimal multiplicity yields reduction number one.

\begin{Theorem}\label{ReductionNumber}
Let $M$ be a  Cohen-Macaulay module of dimension $d$ over a Noetherian local ring $R$ and
let $I$ be an $R$-ideal with $\ell(I, M)=d$. Assume  ${\rm depth}\,(M/IM)\geq {\rm min}\{{\rm dim} (M/IM), 1\}$ and  $I$ satisfies $G_{d}$ and  $AN^-_{d-2}$ on $M$.
 If $I$ has  minimal $j$-multiplicity on $M$
then $r(I,M)=1$.
%and the associated graded module ${\rm gr}_I(M)$ is Cohen-Macaulay.
\end{Theorem}

\demo
By adjoining variables to $R$ and localizing, we may assume  that the residue field is infinite. If ${\rm dim}\,M/IM=0$ then the assertion follows from \cite[Theorem 2.9]{RV}. Now assume ${\rm dim}\,M/IM>0$.
 For general elements $x_1, \ldots, x_d$ in $I$, let  $\overline{M}=M/((x_1, \ldots, x_{d-1})M:_M I^{\infty})$. By Proposition~\ref{lowerbound} and Definition~\ref{j-multiplicity} the $j$-multiplicity can be computed using   $x_1, \ldots, x_d$ thus  $j(I, M)=\lambda(I\overline{M}/I^2\overline{M})$. From  Definition~\ref{j-multiplicity} (see also the proof of Corollary~\ref{rigidity}),\, we obtain $I^2M=x_d IM+ (x_1,\ldots,x_{d-1})M:_{I^2M}I^{\infty}$.
By Lemma  \ref{moduleresidual} (f)
$$(x_1,\ldots,x_{d-1})M:_{I^2M}I^{\infty}=(x_1,\ldots,x_{d-1})IM$$
thus  we  conclude at once that the reduction number of $I$ on $M$ with respect to $(x_1, \ldots, x_d)$ is one. Now  \cite[2.2]{T2} or \cite[8.6.6]{SH} imply  that $r(I,M)=1$.
\QED

\vskip .5cm

Our main application is the case when $M=R$. We obtain that ideals with residual intersection properties and  minimal $j$-multiplicity have Cohen-Macaulay associated graded rings.

\vskip .3cm

\begin{Corollary} \label{idealminimal} Let $R$ be a  Cohen-Macaulay local ring of dimension $d$.  Let $I$ be an $R$-ideal with $\ell(I)=d$. Assume   ${\rm depth}\,(R/I)\geq {\rm min}\{{\rm dim}\,R/I,  1\}$ and $I$ satisfies $G_{d}$ and $AN^-_{d-2}$.
If $I$ has  minimal $j$-multiplicity then the associated graded ring ${\rm gr}_I(R)$ is Cohen-Macaulay.
\end{Corollary}
\demo The assertion follows from  Theorem \ref{ReductionNumber} and \cite[3.1]{JU}.
\QED

\bs

Notice that the `residual intersection assumptions':  $\ell(I)=d$, ${\rm depth}\,(R/I)\geq {\rm min}\{{\rm dim}\,R/I,  1\}$, $I$ satisfies $G_{d}$ and $AN^-_{d-2}$ are all vacuous in the $0$-dimensional case and the condition on the minimal $j$-multiplicity becomes the usual assumption of minimal multiplicity found in the original work of Sally. If the ambient ring is Gorenstein she was able to prove that the associated graded ring is Gorenstein as well. We recover this result in Corollary \ref{Gorenstein}.

We remark that the Artin-Nagata property  $AN^-_{d-2}$ holds if $d \leq g+1$, or if $I$ satisfies $G_d$ and  the sliding depth conditions ${\rm depth}\, R/I^j \geq {\rm
dim}\, R/I - j + 1$ for $1 \leq j \leq d-{\rm ht} \, I + 1$, a linear
weakening of the Cohen-Macaulay property for consecutive powers of
$I$ (see \cite[2.1]{JU}).  The depth inequalities are satisfied by strongly Cohen-Macaulay ideals, i.e., ideals whose Koszul homology modules are Cohen-Macaulay (see \cite[2.10]{U}). Examples of ideals satisfying the latter condition are quite
common -- Cohen-Macaulay almost complete intersections,
Cohen-Macaulay ideals generated by $2+{\rm ht} \, I$ elements,
perfect ideals of codimension two, perfect Gorenstein ideals of
codimension three, and, more generally, any ideal in the linkage class of a complete intersection, namely, licci ideals (see \cite[1.11]{Hu}).

\vskip .3cm

\begin{Corollary}\label{Gorenstein} Let $R$ be a  Gorenstein local ring of dimension $d$.  Let $I$ be an $R$-ideal with $\ell(I)=d$. Assume  ${\rm depth}\, R/I^j \geq {\rm
dim}\, R/I - j + 1$ for $1 \leq j \leq d-{\rm ht} \, I+1$ and $I$ satisfies $G_{d}$ . If $I$ has  minimal $j$-multiplicity then the associated graded ring ${\rm gr}_I(R)$ is Gorenstein.
\end{Corollary}
\demo The assertion follows from  Theorem \ref{ReductionNumber} and \cite[5.3]{JU}.
\QED

\vskip .4cm

\begin{Corollary}\label{SC}  Let $R$ be a  Cohen-Macaulay local ring of dimension $d$.  Let $I$ be a strongly Cohen-Macaulay $R$-ideal with $\ell(I)=d$. If  $I$ has $G_{d}$ and  minimal $j$-multiplicity then the associated graded ring ${\rm gr}_I(R)$ is Cohen-Macaulay. In addition if $R$ is Gorenstein then ${\rm gr}_I(R)$ is Gorensein.
\end{Corollary}
\demo The first assertion follows from  Theorem \ref{ReductionNumber} and \cite[3.2]{JU}, and the second from  \cite[2.10]{U} and Corollary \ref{Gorenstein}.
\QED

\bs

In the next theorem we demonstrate that $r(I,M) \leq 1$ implies the Cohen-Macaulayness of the associated graded module ${\rm gr}_I(M)$.

\begin{Theorem}\label{minimal}
Let $M$ be a  Cohen-Macaulay module of dimension $d$ over a Noetherian local ring $R$ and
let $I$ be an $R$-ideal with $\ell(I, M)=d$. Assume  that $I$ satisfies $G_{d}$ and  $AN^-_{d-2}$ on $M$.
 If $r(I,M) \leq 1$ then  the associated graded module ${\rm gr}_I(M)$ is Cohen-Macaulay.
\end{Theorem}

\demo
By adjoining variables to $R$ and localizing, we may assume that the residue field is infinite.
Set $g= {\rm grade} (I, M)$. Let $x_1,\ldots,x_g$ be general elements in $I$ and $x_1^{*}, \ldots, x_g^{*}$  their initial forms in ${\rm gr}_I(R)$. First we show that
$x_1^{*}, \ldots, x_g^{*}$ form a ${\rm gr}_I(M)$-regular sequence.
By \cite[2.6]{VV}, we only need to show $(x_1,\ldots,x_g)M\cap I^jM=(x_1,\ldots,x_g) I^{j-1}M$ for every $j\geq 1$.
%Since $(x_1,\ldots,x_g)M\cap I^jM\supseteq (x_1,\ldots,x_g) I^{j-1}M$, we just need to show the equality locally at every associated prime of $M/(x_1,\ldots,x_g) I^{j-1}M$. Since

We use induction on $j$ to prove $(x_1,\ldots,x_i)M\cap I^jM=(x_1,\ldots,x_i) I^{j-1}M$ for every $j\geq 1$ and $0\leq i\leq d$. This is clear if $j=1$. So we assume $j\geq 2$ and the equality holds for $j-1$. Set $J=(x_1, \ldots, x_d)$. Since  the reduction number of $I$ on $M$ is one, thus by  \cite[2.2]{T2} or \cite[8.6.6]{SH} we have $I^2M=JIM$. Therefore  $JM\cap I^{j}M=JI^{j-1}M$.
Now we use descending induction on $i$ and assume $(x_1,\ldots,x_{i+1})M\cap I^{j}M=(x_1,\ldots,x_{i+1})I^{j-1}M$. Then
\begin{eqnarray*}
(x_1, \ldots, x_{i})M\cap I^{j}M&=&(x_1, \ldots, x_i)M\cap (x_1,\ldots,x_{i+1})I^{j-1}M \hspace{3.3cm}  \mbox{ by induction on i}\\
&=&(x_1,\ldots,x_i)M\cap ((x_1,\ldots,x_{i})I^{j-1}M+x_{i+1}I^{j-1}M)\\
&=& (x_1,\ldots,x_{i})I^{j-1}M+(x_1,\ldots,x_i)M\cap x_{i+1}I^{j-1}M\\
&=& (x_1,\ldots,x_{i})I^{j-1}M+x_{i+1}[((x_1,\ldots,x_{i})M:_Mx_{i+1})\cap I^{j-1}M]\\
&=&(x_1,\ldots,x_{i})I^{j-1}M+x_{i+1}[(x_1,\ldots,x_{i})M\cap I^{j-1}M]  \ \hspace{.02cm}  \mbox{ by Lemma \ref{moduleresidual} (a)  and (e)}\\
&= &(x_1,\ldots,x_{i})I^{j-1}M+x_{i+1}(x_1,\ldots,x_{i}) I^{j-2}M \hspace{2.3cm}  \mbox{ by induction on j}\\
&\subseteq &(x_1,\ldots,x_i)I^{j-1}M.
\end{eqnarray*}

Set  $\delta(I,M)=d-g$. Now we prove that the associated graded module ${\rm gr}_I(M)$ is Cohen-Macaulay by induction on $\delta$.  If  $\delta=0$, the assertion follows because $x_1^{*}, \ldots, x_g^{*}$ form a ${\rm gr}_I(M)$-regular sequence.  Thus we may assume  $\delta(I,M)\geq 1$ and the theorem holds for smaller values of $\delta(I,M)$. In particular, $d\geq g+1$. Again since $x_1^{*}, \ldots, x_g^{*}$ form a ${\rm gr}_I(M)$-regular sequence, we may factor out $x_1, \ldots, x_g$ to assume $g=0$. Now $d=\delta(I,M) \geq 1$.
Set  $H_0=0_M:_M I$ and $\overline{M}=M/H_0$. Then $\overline{M}$ is Cohen-Macaulay since $I$ satisfies $AN^-_{d}$ on $M$ (see Lemma \ref{Moduleresidual}).
By Lemma \ref{moduleresidual} (e), (a) and (g) and Lemma \ref{Moduleresidual}, $IM\cap H_0=0,\, {\rm grade}\,(I, \overline{M})\geq 1$, $I$ still satisfies  $G_{d}$ and $AN^-_{d}$ on $\overline{M}$.
Furthermore ${\rm dim}\, M={\rm dim}\,\overline{M}=d$ and $IM\cap H_0=0$ implies $\ell(I, \overline{M})=\ell(I, M)=d$.  Again by $IM\cap H_0=0$, there is a graded exact sequence
\begin{equation}\label{eq11}
0\rightarrow H_0\rightarrow {\rm gr}_I(M)\rightarrow {\rm gr}_{I}(\overline{M})\rightarrow 0.
\end{equation}
%Notice that   ${\rm depth}\,H_0=d$ since ${\rm depth}\,\overline{M}={\rm depth} \, M=d$. By the above exact sequence in degree zero
%\begin{eqnarray*}{\rm depth}\,(\overline{M}/I\overline{M})&\geq& {\rm min}\{{\rm depth}\,H_0-1,\,{\rm depth}\,M/IM\}\\
 %&\geq& {\rm min}\{d-{\rm grade}\,(I, \overline{M}),\,1\}\\
 %&=& {\rm min}\{{\rm dim} (\overline{M}/I\overline{M}),\,1\}.
 %\end{eqnarray*}
%  Observe that $x_1, \ldots, x_d$ are still general elements of $I$ on $\overline{M}$ thus we can compute the  $j$-multiplicity of $I$ on   $\overline{M}$ using them. Once again $IM \cap H_0=0$ yields
  %$$\lambda (I^2\overline{M}/[ (x_1, \ldots, x_{d-1})\overline{M}:_{I^2\overline{M}}I^{\infty}+x_d I\overline{M} ])=\lambda (I^2M/[ (x_1,\ldots,x_{d-1})M:_{I^2M}I^{\infty}+x_d IM])= 0.$$ In particular, $I$ has minimal $j$-multiplicity on $\overline{M}$.
Notice that $r(I, \overline{M}) \leq 1$.  Since $\delta(I,\overline{M})=d-{\rm grade}\,(I, \overline{M})<d=\delta(I,M)$, by induction hypothesis  ${\rm depth}({\rm gr}_{I}(\overline{M}))\geq d$.  Observe ${\rm depth}(H_0)\geq d$ since $\overline{M}$ is Cohen-Macaulay. The Cohen-Macaulyness of ${\rm gr}_I(M)$ follows at once by the short  exact sequence (\ref{eq11}). \QED
\ms

Now we are ready to prove the main theorem.

\begin{Theorem}\label{MAIN}
Let $M$ be a  Cohen-Macaulay module of dimension $d$ over a Noetherian local ring $R$ and
let $I$ be an $R$-ideal with $\ell(I, M)=d$. Assume  ${\rm depth}\,(M/IM)\geq {\rm min}\{{\rm dim} (M/IM), 1\}$ and  $I$ satisfies $G_{d}$ and  $AN^-_{d-2}$ on $M$.
 If $I$ has  minimal $j$-multiplicity on $M$
then  the associated graded module ${\rm gr}_I(M)$ is Cohen-Macaulay.
\end{Theorem}
\demo
By Theorem \ref{ReductionNumber}, the reduction number of $I$ on $M$ is one. Now the assertion follows from Theorem \ref{minimal}.
\QED

\section{Almost minimal $j$-multiplicity}

\ms

We start by giving the definition of almost minimal $j$-multiplicity, which is the analogue of almost minimal multiplicity  \cite{RV}.

\begin{Definition}\label{almost j-multiplicity}
{\rm  Let $M$ be a finite module of dimension $d$ over a Noetherian local ring $R$ and
  $I$ an $R$-ideal with analytic spread $\ell(I, M)=d$. We say that $I$ has {\it almost minimal $j$-multiplicity on  $M$} if  $j(I, M)= \lambda(I\overline{M}/I^2\overline{M})+1$, where $\overline{M}= M/((x_1,\ldots,x_{d-1})M:_M I^{\infty})$ and $x_1, \ldots, x_{d-1}$ are general in $I$.}
\end{Definition}

\ms
Notice that by Lemma \ref{general} the definition of almost minimal $j$-multiplicity is independent on the general sequence chosen in $I$.

\begin{Remark}\label{AMM}{\rm  If $I$ has almost minimal $j$-multiplicity on  $M$ then
$$
\lambda(I^2M/[x_d IM+ (x_1,\ldots,x_{d-1})M:_{I^2M}I^{\infty}])=1.
$$
for any general sequence $x_1, \ldots, x_d$ in $I$. }
\end{Remark}

\ms

%\begin{Lemma}
%Let $R$ be a $1$-dimensional Noetherian local ring  and $I$  an $R$-ideal. Let $M$  be a finite $R$-module with $(0_M:_M I^{\infty})\cap IM=0_M$. Let $G={\rm gr}_I(R)$ and $T={\rm gr}_I(M)$. Then for every $j\geq 1$, $H(j)=\lambda(I^jM/I^{j+1}M)<\infty$. Let $e_0=H(j)$ for large $j$. Define for every $j\geq 1$, $\nu_j=e_0-H(j)$. Let $x\in I\setminus I^2$ such that the initial form $x^*$ is filter-regular for $T$ with respect to $G_+$, then
%for every $j\geq 1$,
%$$
%\nu_j=\lambda(I^{j+1}M/xI^jM).
%$$
%\end{Lemma}

%\demo
%For sufficiently large $j$, there exist exact sequences:
%$$
%0\rightarrow I^{j+1}M:_{IM}x/I^jM \rightarrow IM/I^jM \rightarrow IM/I^{j+1}M\rightarrow IM/(I^{j+1}M+xIM) \rightarrow 0, $$
%and
%$$
%0\rightarrow (0_M:_{IM}x)/(0_M:_{I^jM}x) \rightarrow IM/I^jM \rightarrow xIM/xI^jM \rightarrow 0.
%$$
%Then we have $e_0=\lambda(IM/xIM)=\lambda(IM/xI^jM)-\lambda(xIM/xI^jM)=\lambda(IM/xI^jM)-\lambda(IM/I^jM)$ and $\nu_j=e_0-H(j)=\lambda(IM/xI^jM)-\lambda(IM/I^jM)-\lambda(I^jM/I^{j+1}M)=\lambda(I^{j+1}M/xI^jM)$.
%\QED

%\ms

%\begin{Lemma}\label{Hilbertcoefficient}
%Let $R$ be a $1$-dimensional Noetherian local ring  and $I$  an $R$-ideal. Let $M$  be a finite $R$-module with $(0_M:_M I^{\infty})\cap IM=0_M$. Let $G={\rm gr}_I(R)$ and $T={\rm gr}_I(M)$.
%Let $P(z)=\sum_{j=1}^{\infty}H(j)z^j=h(z)/(1-z)$, where $h(z)=h_1z+\cdots h_sz^s$.
%Then
%$$e_1=\sum_{j=2}^s (j-1)h_j=\sum_{j=1}^{s-1}\nu_j=\sum_{j= 1}^{s-1}\lambda(I^{j+1}M/xI^jM).$$
%\end{Lemma}

%\ms

\begin{Notation and Discussion}\label{RR}Let $M$ be a  finite module over a Noetherian local ring $R$ and  $I$  an $R$-ideal. For every $j\geq 1$, let $\widetilde{I^j}M=\cup_{t\geq 1}(I^{j+t}M:_{M}I^t)$ be the {\it Ratliff-Rush} filtration of $I$ on $M$
$($see  \cite{RR, PZ, RV}$)$.
If ${\rm depth}_I\,M > 0$, by \cite[Lemma 3.1]{RV}, there exists an integer $n_0$ such that  $\widetilde{I^j}M= I^jM$ for every $j \geq n_0$. In particular, for every  reduction $J$ of $I$ on $M$ and every $j \geq n_0$, we obtain $\widetilde{I^{j+1}}M= J\widetilde{I^j}M+ I^{j+1}M$.
Thus the module $N:= \oplus_{j \geq 0} (\widetilde{I^{j+1}}M/ J\widetilde{I^j}M+ I^{j+1}M)$ has finitely many non-zero components. Since each component is a finitely generated $R$-module, $N$ itself is finitely generated as an $R$-module; we denote with $q=\mu(N)$ its minimal number of generators.
\end{Notation and Discussion}
\ms

The next result is the key step in relating the reduction number of an $I$-adic filtration to invariants of the Ratliff-Rush filtration of $I$ on $M$. The idea originated in the work of  Rossi and  Valla (see \cite{RV1, R1,  RV}). For clarity of
exposition we state here the version for modules and ideals that are not necessarily ${\mathfrak m}$-primary.

\begin{Theorem}
Use the notation of \ref{RR} and assume   ${\rm depth}_I\,M > 0$. Let $J$ be an ideal generated by $d$ general elements in $I$.  Then $$I^q \subseteq J I^{q-1} + (I^{q+j}M :_R \widetilde{I^j}M)$$ for every positive integer $j$.
\end{Theorem}
\demo
The proof is the same as the proof of  \cite[Theorem 4.1]{RV}.
\QED

\ms
\begin{Corollary}\label{red}
Use the notation of \ref{RR} and assume   ${\rm depth}_I\,M > 0$. Let $J$ be an ideal generated by $d$ general elements in $I$. Then  $r(I,M) \leq t+q$, where $t={\rm min} \{j \mid I^{j+1}M \subseteq J \widetilde{I^j}M \}$.
\end{Corollary}
\demo
The proof is the same as the proof of \cite[Corollaries 4.1 and 4.2]{RV}.
\QED
\ms

\begin{Lemma}\label{length}
Let $M$ be a  Cohen-Macaulay module of dimension $d$ over a Noetherian local ring $R$.
Let $I$ be an $R$-ideal with $\ell(I, M)=d$. Assume  that  $I$ satisfies $G_{d}$ and $AN^-_{d-2}$ on $M$ and let $J$ be an ideal generated by $d$ general elements  in $I$.  Then the lengths  $\lambda( I^j IM/J I^{j-1} IM) $,  $\lambda( \widetilde{I^j} IM/ J\widetilde{ I^{j-1}} IM) $
and  $\lambda( \widetilde{I^j} IM/ I^j IM) $  are finite for all $j \geq 1$. \end{Lemma}
\demo
Clearly $\lambda(  IM/JM) < \infty $ and $\lambda( I^j IM/J I^{j-1} IM) < \infty$   for every $j \geq 1$,   since $IM$ and $JM$ are the same on the punctured spectrum.
To show the remaining assertions, observe that $J$ has analytic spread $d$ and satisfies $G_{d}$ and $AN^-_{d-2}$ on $M$ as well, for instance by Lemma ~\ref{generalresidual} and \cite[3.1]{U}. Since $r(J,M)=0$, by Theorem \ref{minimal} the associated graded module ${\rm gr}_J (M)$ is  Cohen-Macaulay.
In particular, $J$ is Ratliff-Rush closed on $M$, i.e., $\widetilde{J^j} M=J^jM$ for all $j \geq 1$.
Thus on the punctured spectrum $I$ is Ratliff-Rush closed on $M$ as well, in particular $\lambda( \widetilde{I^j} M/ J^j M) < \infty$ for every $j\geq 1$.
\QED

\ms
\begin{Theorem}\label{dim2}
Let $M$ be a  Cohen-Macaulay module of dimension $2$ over a Noetherian local ring $R$ with infinite residue field.
Let $I$ be an $R$-ideal with $\ell(I, M)=2$. Assume   ${\rm depth}\,(M/IM)\geq {\rm min}\{{\rm dim} (M/IM), 1\}$ and  $I$ satisfies $G_{2}$ and  $AN^-_{0}$ on $M$. Let $x_1$ be a general element in $I$.  If $I$ has almost minimal $j$-multiplicity on  $M$ then
\begin{itemize}

\item[(a)] $x_1^*$ is regular on ${\rm gr}_I(M)_+$ $;$
\item[(b)]${\rm depth}\,( {\rm gr}_I(M) )\geq 1$ .
\end{itemize}
\end{Theorem}

\demo
 If ${\rm dim}\,M/IM=0$ then both claims follow from \cite[Theorem 4.4]{RV}. Thus we may assume ${\rm depth}\,(M/IM)>0$.

Since $I$ has almost minimal $j$-multiplicity on  $M$, by Remark \ref{AMM}, for general elements
 $x_1,\, x_2$ in $I$,   $\lambda (I^2M/[x_2 IM+ x_1M:_{I^2M}I^{\infty}])= 1$. Set  $J =(x_1, x_2)$. By Lemma \ref{moduleresidual} (f), we have  $x_1M:_{I^2M}I^{\infty}=x_1 IM$, thus $\lambda(I^2M/JIM)=1$.
Hence $I^2M=ab+JIM$ for some $a\in I, b\in IM$ with $ab\not\in JIM$. For $j\geq 2$, the multiplication by $a$ gives a surjective map from $I^jM/JI^{j-1}M$ to $I^{j+1}M/JI^{j}M$. Thus $\lambda (I^jM/JI^{j-1}M)\leq 1$ for every $j\geq 2$.

%First if  ${\rm grade}(I, M)=0$, let $H_0=0_M:_M I$ and $\overline{M}=M/H_0$. As in the proof of Theorem \ref{minimal},  ${\rm grade}(I, \overline{M})\geq 1$, $\overline{M}$ is a  Cohen-Macaulay module of dimension $2$, $\ell(I, \overline{M})=2$,   and  $I$ satisfies $G_{2}$ and  $AN^-_{0}$ on $\overline{M}$. Since  $IM\cap H_0=0$, ${\rm gr}_I(M)_+={\rm gr}_I(\overline{M})_+$. By the exact sequence
%$$
%0\rightarrow H_0\rightarrow {\rm gr}_I(M)\rightarrow {\rm gr}_{I}(\overline{M})\rightarrow 0,
%$$
%we have
%\begin{eqnarray*}
%{\rm depth}\,(\overline{M}/I\overline{M})& \geq & {\rm min}\{{\rm depth}\,H_0-1,\,{\rm depth}\,M/IM\}\\
 %& \geq & {\rm min}\{2-{\rm grade}\,(I, \overline{M}),\,1\}\\
 %&= &{\rm min}\{{\rm dim} (\overline{M}/I\overline{M}),\,1\}.
 %\end{eqnarray*}

%Again as  in the proof of Theorem \ref{minimal}, the assertions  will follow immediately once we proved that (a) and (b) hold for $\overline{M}$. Thus we reduce to the case where ${\rm grade}(I, M)>0$. Now if ${\rm dim}\,M/IM=0$ then the results follow from \cite[Theorem 4.4]{RV}. So we may assume ${\rm grade}(I, M)=1$, hence ${\rm depth}\,(M/IM)>0$.

Notice that $x_1$ is regular on $IM$ since $(0:_M x_1)  \cap IM =0$ (Lemma \ref{moduleresidual} (e)).  Thus to prove that  $x_1^*$ is regular on ${\rm gr}_I(M)_+={\rm gr}_I(IM)$ we only need to show $x_1 IM\cap I^jIM=x_1 I^{j-1}IM$
for every $j\geq 1$ by \cite[2.6]{VV} (see also \cite[Lemma 1.1]{RV}).
This is clear if $j=1$; hence we can  assume $j\geq 2$. Let $^{^{\tratto}}$ denote images
in $\overline{IM}=IM/x_1M$ and set  $s=r_{J}(I, \overline{IM})$. We claim that it is enough to show   $r_{J}(I, IM)=s$. Indeed,
if $1\leq j\leq s$ then $JI^{j-1}IM+(x_1M\cap I^{j}IM)=JI^{j-1}IM$. This follows from the following easy inequality of lengths
\vskip -.4cm
\begin{eqnarray*}
0&<&\lambda (I^j IM/ JI^{j-1}IM+(x_1M \cap I^{j} IM))\\
&=&\lambda(I^j IM/ JI^{j-1}IM) - \lambda(JI^{j-1}IM+(x_1M \cap I^{j}IM) / JI^{j-1}IM)\\
&=& 1 - \lambda(JI^{j-1}IM+(x_1M \cap I^{j}IM)/JI^{j-1}IM).
\end{eqnarray*}

\vskip .2cm

\noindent
On the other hand, if $j\geq s+1= r_{J}(I, IM)+1$, then  $I^jIM=JI^{j-1}IM$. Thus  we have for all $j\geq 1$

\vskip -.2cm
\begin{equation}\label{eq1}
x_1IM\cap I^{j}IM=x_1IM\cap JI^{j-1}IM.
\end{equation}

\vskip .2cm
\noindent
Now we  proceed by induction on $j \geq 2$ as in the proof of  Theorem \ref{minimal}:
\begin{eqnarray*}
x_1IM\cap I^{j}IM&=&x_1IM\cap JI^{j-1}IM
\hspace{4cm}  {\rm by} \   (\ref{eq1}) \\
&=&x_1 IM \cap (x_1 I^{j-1} IM + x_{2}I^{j-1} IM)\\
&=& x_1 I^{j-1} IM + (x_1 IM \cap x_{2}I^{j-1} IM)\\
&=&  x_1 I^{j-1} IM + x_{2}[(x_1 IM:_{IM} x_{2}) \cap I^{j-1} IM] \\
&=& x_1 I^{j-1} IM + x_{2}[x_1 IM\cap I^{j-1}IM].
\end{eqnarray*}

The last equality follows because $x_1 IM\cap I^{j-1}IM \subseteq ( x_1I M:_{IM} x_2)  \cap I^{j-1} IM \subseteq (x_1M:_{I^2M}x_2)  \cap I^{j-1} IM = x_1 IM\cap I^{j-1}IM$ because $j \geq 2$ and by Lemma \ref{moduleresidual} (f).
Now using induction on $j$, we obtain
$$\hspace{.5cm} =x_1 I^{j-1} IM + x_{2}[x_1I^{j-2}IM] = x_1 I^{j-1}IM.$$

To complete the proof of (a), we still need to to show that $r_{J}(I, IM)=s$. For this purpose we will use the Ratliff-Rush filtration $\widetilde{I^j} IM$ as it is done for ideals of definition
(see \cite[Theorem 4.2]{RV}). As noticed earlier  $x_1$ is regular on $IM$. Thus, for instance by  \cite[Lemma 3.1]{RV},
%${\rm depth}_I\, IM \geq 1$. The latter condition implies that
there exists an integer $n_0$ such that
$I^jIM=\widetilde{I^j} IM$ for $j\geq n_0$, and
\begin{equation}\label{RRF}
\widetilde{I^{j+1}} IM  :_{IM} x_1=\widetilde{I^{j}}IM \quad \mbox{for every } j\geq 0.
\end{equation}

As before, let $\overline{IM}=IM/x_1M$ and $^{^{\tratto}}$   denote images in $\overline{IM}$.  There are two filtrations:
$$
\overline{\mathbb{M}}: \overline{IM}\supseteq  \overline{I^2M}=I\,\overline{IM}\supseteq\ldots\supseteq I^{j-1} \overline{IM}\supseteq \ldots
$$
and
$$
\overline{\mathbb{N}}:\overline{IM} \supseteq \overline{\widetilde{I} IM}= \widetilde{I} \   \overline{IM}\supseteq\ldots\supseteq \widetilde{I^{j-1}} \overline{IM}\supseteq \ldots
$$
Notice that $\overline{\mathbb{M}}$ is an  $I$-adic filtration and $\overline{\mathbb{N}}$ is a good $I$-filtration on $\overline{IM}$ (see \cite[Page 9]{RV} for the definition of  good filtration).
Furthermore, $I$ is an ideal of definition for $\overline{IM}$, i.e., $\lambda_{R}\,(\overline{IM}/  I\overline{IM}) < \infty $. Indeed, $(x_1M :_{M} x_2) \cap IM =x_1  M$  (see Lemma \ref{moduleresidual} (e)) which forces
$x_2 \in I$
to be regular on $\overline{IM}$, in turns this  yields $\lambda_{R}(\overline{IM}/
I\overline{IM})  \leq \lambda_{R}(\overline{IM}/ x_2 \overline{IM}) < \infty. $
% Indeed, since  $\overline{R}$  is a $1$-dimensional Noetherian local ring
%and  $\overline{I}$ satisfies  condition $G_1$ (see Lemma \ref{}), we have $(0_{\overline{M}}:_{\overline{M}}\overline{I}^{\infty})\cap \overline{I}\overline{M}=0_{\overline{M}}$. Hence, there exist an $\xi \in (notice that we can take $\xi=x_2$
%because $0_{\overline{M}}:_{\overline{M}} \overline{I}^{\infty} = 0_{\overline{M}}:_{\overline{M}} \, x_2$ by Lemma \ref{eq}).
 Thus we are  in the context of the filtrations treated in \cite{RV}.
Since $ I^{j-1} \, \overline{IM}=\widetilde{I^{j-1}}\overline{IM}$ for $j\geq n_0$, the associated graded modules  ${\rm gr}_{\overline{\mathbb{M}}}(\overline{IM})$ and ${\rm gr}_{\overline{\mathbb{N}}}(\overline{IM})$
have the same Hilbert coefficients $e_0$ and $e_1$. Again, because there exists an element in  $I$
which is regular on $\overline{IM}$, by \cite[Lemmas 2.1 and 2.2]{RV} we have

\begin{equation} \label{eq3} s=\sum_{j\geq 0} \lambda(I^{j+1}\overline{IM}/ x_2 I^{j}\overline{IM}) =e_1(\overline{\mathbb {M}})=e_1(\overline{\mathbb {N}}) =
\sum_{j\geq 0} \lambda(\widetilde{I^{j+1}}\overline{ IM}/ x_2\widetilde{I^{j}} \overline{ IM}). \end{equation}

 Observe that the first equality holds  because
 $ \lambda(I^{j+1}\overline{IM}/ x_2 I^{j}\overline{IM}) =1$ for all $0 \leq j \leq s-1$ and  $  \lambda(I^{j+1}\overline{IM}/ x_2 I^{j}\overline{IM}) =0$ for $j \geq s= r_{J}(I, \overline{IM})$.

 \medskip

We  prove that $ \lambda(\widetilde{I^{j+1}}\overline{ IM}/ x_2\widetilde{I^{j}} \overline{ IM}) =  \lambda(\widetilde{I^{j+1}} IM/ J \widetilde{I^j} IM)$ for every $j\geq 0$. Since
$$\widetilde{I^{j+1}}\overline{ IM}/ x_2\widetilde{I^{j}} \overline{ IM}\cong \widetilde{I^{j+1}} IM/(x_1M\cap \widetilde{I^{j+1}}IM+ x_2\widetilde{I^{j}} IM),$$ we just need to show $x_1M\cap \widetilde{I^{j+1}}IM=x_1\widetilde{I^{j}} IM$. We first prove $x_1M\cap \widetilde{I}IM=x_1IM$. Since $x_1M\cap \widetilde{I}IM\supseteq x_1IM$, it suffices to show the equality locally at every associated prime ideal of $M/x_1IM$. By Lemma \ref{moduleresidual} (d), every $\p\in {\rm Ass}(M/x_1IM)$ is not maximal. By the proof of Lemma \ref{length}, $x_1M_\p=\widetilde{I}M_\p=I M_\p$ thus $x_1 M_\p\cap \widetilde{I} IM_\p= \widetilde{I} IM_\p=x_1I M_\p$. Therefore $x_1M\cap \widetilde{I}IM=x_1IM$. Now for any $j\geq 1$, $x_1M\cap \widetilde{I^{j+1}}IM=x_1IM \cap \widetilde{I^{j+1}}IM=x_1(\widetilde{I^{j+1}}IM :_{IM} x_1)=x_1\widetilde{I^{j}} IM$, where the last equality holds by (\ref{RRF}).

Now  (\ref{eq3})  gives us

\begin{equation}\label{eq4}
\sum_{j \ge 0}\lambda(\widetilde{I^{j+1}} IM/ J \widetilde{I^j} IM) =s,
\end{equation}
and
\begin{equation*}p=  {\rm inf} \{j  \mid J\widetilde{I^{j}} IM=\widetilde{I^{j+1}}  IM \}\leq s.   \end{equation*}

Let $t=  {\rm inf} \{j  \mid I^{j+1} IM \subseteq J \widetilde{I^{j}}  IM \}$. Observe that $t \leq p \leq s$ because $I^{p+1} IM \subseteq \widetilde{I^{p+1}}  IM =J\widetilde{I^{p}} IM$. Let $l$ be a positive integer such that for all $0 \leq j \leq l$ we have
$I^{j+1} IM \cap J IM = J I ^{j} IM$.  If $t \leq l$, then  $r_{J}(I, IM) \leq t  \leq s$ and we are done. So we can assume that $t > l$ and we have:

\vskip -.4cm
$$ l< t \leq p \leq s.
$$

By Corollary \ref{red}, the reduction number is bounded above by $t+q$ where

\begin{equation}\label{eq5}
q \leq   \sum_{j \geq 0} \lambda(\widetilde{I^{j+1}}IM/ J\widetilde{I^j}IM+ I^{j+1}IM).
\end{equation}

To prove that  $r_{J}(I, IM)=s$, it will be enough to show that  $ \sum_{j \geq 0} \lambda(\widetilde{I^{j+1}}IM/ J\widetilde{I^j}IM+ I^{j+1}IM) \leq s-t$. Observe that for $0 \leq j \leq l$, we can relate  $\lambda(\widetilde{I^{j+1}}IM/ J\widetilde{I^j}IM+ I^{j+1}IM)$ to the difference of the length of the factors of the filtrations ${\mathbb N}$ and ${\mathbb M}$. Indeed, for $0 \leq j \leq l$ we have $J \widetilde{I^j} IM \cap I^{j+1} IM= J I ^{j} IM$ and therefore we obtain the following family of short exact sequences:
$$
0 \longrightarrow J \widetilde{I^j}IM/ J I^j IM \longrightarrow \widetilde{I^{j+1}}IM/I^{j+1} IM \longrightarrow \widetilde{I^{j+1}}IM/ J\widetilde{I^j}IM+ I^{j+1}IM \longrightarrow 0,
$$
from which we obtain:

\begin{eqnarray*}
\lambda(\widetilde{I^{j+1}}IM/ J\widetilde{I^j}IM+ I^{j+1}IM) &=&  \lambda(\widetilde{I^{j+1}}IM/I^{j+1} IM) - \lambda(J  \widetilde{I^j}IM/ J I^j IM) \\
&=&\lambda( \widetilde{I^{j+1}} IM/ J \widetilde{I^j}IM) -  \lambda(I^{j+1} IM/J I^{j} IM) \quad \quad  \quad  \mbox{for} \  0 \leq j \leq l.
\end{eqnarray*}

For this we conclude that
\begin{eqnarray}
\lambda(\widetilde{I^{j+1}}IM/ J\widetilde{I^j}IM+ I^{j+1}IM) &=& \lambda( \widetilde{I^{j+1}} IM/ J \widetilde{I^j}IM) - 1 \quad  0 \leq j \leq l, \\
\lambda(\widetilde{I^{j+1}}IM/ J\widetilde{I^j}IM+ I^{j+1}IM) &\le& \lambda( \widetilde{I^{j+1}} IM/ J \widetilde{I^j}IM) -1 \quad  l+1 \le j \le t-1,\\
\lambda(\widetilde{I^{j+1}}IM/ J\widetilde{I^j}IM+ I^{j+1}IM)  &=& \lambda( \widetilde{I^{j+1}} IM/ J \widetilde{I^j}IM) \quad \quad  \quad  j \ge t.
\end{eqnarray}

\noindent
Now by means of (\ref{eq4}), (\ref{eq5}),  (6), (7) and (8) we obtain
$$
q \leq   \sum_{j \geq 0} \lambda(\widetilde{I^{j+1}}IM/ J\widetilde{I^j}IM+ I^{j+1}IM) \le \sum_{j \ge 0}\lambda(\widetilde{I^{j+1}} IM/ J \widetilde{I^j} IM)  -t =s -t.
$$
\noindent
This concludes the proof of (a) since $r_{J}(I, IM)\le t+q \le t+ s -t=s$ (see Corollary \ref{red}).

\medskip

% If ${\rm grade}\,(I, M)=2$, then $I$ is an ideal of definition for $M$ and both assertions follow from \cite[4.2.1]{RV}.

% We need to show that  ${\rm depth}\,( {\rm gr}_I(M) )\geq 1$.
Finally part (b) follows from  (a). Indeed,  by assumption  ${\rm depth}\,(M/IM)>0$ and by the exact sequence:
$$
0\rightarrow M/IM \rightarrow {\rm gr}_I(M) \rightarrow {\rm gr}_I(M)_+ \rightarrow 0,
$$
we conclude ${\rm depth} ({\rm gr}_I(M))\geq {\rm min} \{{\rm depth} \, M/IM , {\rm depth} ({\rm gr}_I(M)_+) \} \geq 1$.

\QED

\ms
\begin{Theorem}
Let $M$ be a  Cohen-Macaulay module of dimension $d$ over a Noetherian local ring $R$.
Let $I$ be an $R$-ideal with $\ell(I, M)=d$. Assume   ${\rm depth}\,(M/IM)\geq {\rm min}\{{\rm dim} (M/IM), 1\}$ and  $I$ satisfies $G_{d}$ and  $AN^-_{d-2}$ on $M$.
 If $I$ has almost minimal $j$-multiplicity on $M$
then ${\rm depth}\,( {\rm gr}_I(M) )\geq d-1.$
%\begin{itemize}
%\item[(a)] $x_1^*, \ldots, x_{d-1}^*$ is a regular sequence on ${\rm gr}_I(M)_+$;
%\item[(b)]${\rm depth}\,( {\rm gr}_I(M) )\geq d-1$.
%\end{itemize}
\end{Theorem}

\demo As in the proof of Theorem \ref{ReductionNumber},  we may
assume  that the residue field of $R$ is infinite.   We prove the
theorem by induction on $d$. The case  $d=2$ being proven in Theorem
\ref{dim2}. Let $d\geq 3$ and assume the theorem holds for $d-1$. We
first reduce to the case  ${\rm grade}\, (I, M)  \geq1$. If ${\rm
grade}\,(I, M)=0$,  let $H_0=0:_M I$. As in the proof of Theorem
\ref{minimal}, all assumptions still hold for the module $M/H_0$.
Furthermore  $IM/H_0=IM/(H_0 \cap IM)=IM$, ${\rm grade} \, (I,
M/H_0) \geq 1$ and again as in the proof of  Theorem \ref{minimal},
${\rm depth}\,( {\rm gr}_I(M) )\geq {\rm depth}\,( {\rm
gr}_{I}(M/H_0) )$. So we are reduced to the case where the ideal $I$
has at least one  $M$-regular element.  Thus if $x_1$ is a general
element in $I$ then $x_1$ is regular on $M$.

If  ${\rm dim}\,M/IM=0$ then the assertion  follows from
\cite[Theorem 4.4]{RV}. Thus we may assume ${\rm dim}\,M/IM>0$. Let
$^{^{\tratto}}$  denote images in $\overline{M}=M/x_1M$.   Observe
that $\overline{M}$ is a Cohen-Macaulay module of dimension $d-1$
and $\ell(I, \overline{M})=d-1$. Also $I$ satisfies $G_{d-1}$ and
$AN^-_{d-3}$ on $\overline{M}$ by Lemma \ref{moduleresidual}.
Furthermore, observe  $\overline{M}/I\overline{M} \cong M/IM $  thus
${\rm depth}\,(\overline{M}/\overline{IM})={\rm depth}\,(M/IM)\geq
{\rm min}\{{\rm dim} \, M/IM, 1\}= \{{\rm dim} \,
\overline{M}/I\overline{M}, 1\}$. Clearly  $I$ has almost minimal
$j$-multiplicity on $\overline{M}$. By induction hypothesis, ${\rm
depth}\,( {\rm gr}_{I}(\overline{M}) )\geq d-2.$

 Now we prove that $x_1^*$ is regular on ${\rm gr}_I(M)$. Since $x_1$ is regular on $M$, by \cite[2.6]{VV} (see also \cite[Lemma 1.1]{RV}), the claim follows if the intersections $x_1M\cap I^{j}M=x_1I^{j-1}M$  hold for every $j\geq 1$. This is clear if $j=1$. If $j=2$, since $x_1IM\subseteq x_1M\cap I^2M$, it suffices to show the equality locally at every prime ideal $\p\in {\rm Ass}(M/x_1IM)$. By  Lemma \ref{moduleresidual} (d), ${\rm depth}\,(M/x_1IM)\geq 1.$ Thus for every prime ideal  $\p\in {\rm Ass}(M/x_1IM)$, $\p$ is not the maximal ideal of $R$ and hence either $IM_{\p}=M_{\p}$ or $IM_{\p}=(x_1, \ldots, x_{d-1})M_{\p}$. Therefore $(x_1,\ldots,x_{d-1})M_{\p}\cap I^2M_{\p}=(x_1, \ldots, x_{d-1})IM_{\p}$. We use descending induction on $i$ to prove $(x_1,\ldots,x_{i})M_{\p}\cap I^2M_{\p}=(x_1, \ldots, x_{i})IM_{\p}$ for every $1\leq i\leq d-1$. Assume $(x_1,\ldots,x_{i+1})M_{\p}\cap I^2M_{\p}=(x_1,\ldots,x_{i+1})IM_{\p}$. Then
\begin{eqnarray*}
(x_1, \ldots, x_{i})M_{\p}\cap I^2M_{\p}&=&(x_1, \ldots, x_i)M_{\p}\cap (x_1,\ldots,x_{i+1})IM_{\p}\\
&=&(x_1,\ldots,x_i)M_{\p}\cap ((x_1,\ldots,x_{i})IM_{\p}+x_{i+1}IM_{\p})\\
&=& (x_1,\ldots,x_{i})IM_{\p}+(x_1,\ldots,x_i)M_{\p}\cap x_{i+1}IM_{\p}\\
&=& (x_1,\ldots,x_{i})IM_{\p}+x_{i+1}[((x_1,\ldots,x_{i})M_{\p}:_{M_{\p}}x_{i+1})\cap IM_{\p}]\\
%&=&(x_1,\ldots,x_{i})IM_{\p}+x_{i+1}[(x_1,\ldots,x_{i})M_{\p}\cap IM_{\p}] \\
&=& (x_1,\ldots,x_{i})IM_{\p}+x_{i+1}(x_1,\ldots,x_{i}) M_{\p}\subseteq (x_1,\ldots,x_i)IM_{\p}.
\end{eqnarray*}
When $j\geq 3$, we have $x_1M\cap I^{j}M=x_1M\cap I^2M \cap I^{j}M=x_1IM\cap I^{j-1}IM=x_1I^{j-2}IM=x_1I^{j-1}M$ since $x_1^*$ is regular on ${\rm gr}_I(M)_+$ by Theorem \ref{dim2}.  Finally since ${\rm depth}({\rm gr}_{I}(\bar{M}))\geq d-2$ and $x_1^*$ is regular on ${\rm gr}_I(M)$, we have ${\rm depth}({\rm gr}_{I}(M))\geq d-1$.

\QED

\ms

Again our main application is the case when $M=R$. We obtain that ideals with residual intersection properties and almost minimal $j$-multiplicity have associated graded rings almost Cohen-Macaulay.

\begin{Corollary} \label{idealalmostminimal} Let $R$ be a  Cohen-Macaulay local ring of dimension $d$.  Let $I$ be an $R$-ideal with $\ell(I)=d$. Assume   ${\rm depth}\,(R/I)\geq {\rm min}\{{\rm dim}\,R/I,  1\}$ and $I$ satisfies $G_{d}$ and $AN^-_{d-2}$. If $I$ has almost minimal $j$-multiplicity then ${\rm depth}\,( {\rm gr}_I(R) )\geq d-1$.
\end{Corollary}

\ms

As noticed in Section 3 before Corollary \ref{Gorenstein},  the assumptions on Corollary \ref{idealalmostminimal} are all vacuous in the $0$-dimensional case and the condition on the almost minimal $j$-multiplicity becomes the usual assumption of almost minimal multiplicity found in Sally's conjecture  and in the work of Rossi and Valla and Wang (see \cite{S3}, \cite{RV1}, \cite{W}, \cite{H}, \cite{CPV}, \cite{R1}, \cite{E}). Thus Corollary \ref{idealalmostminimal} can be viewed as a positive answer to Sally's conjecture for arbitrary ideals.

Again the conclusion of Corollary \ref{idealalmostminimal} will hold true for $I$  which is generically a complete intersection with ${\rm ht} \, I =d-1$ or for  strongly Cohen-Macaulay ideals satisfying $G_d$, in particular, for  Cohen-Macaulay almost complete intersections,
Cohen-Macaulay ideals generated by $2+{\rm ht} \, I$ elements,
perfect ideals of codimension two, perfect Gorenstein ideals of
codimension three, and, more generally licci ideals.

We will finish our paper by the following examples.

\ms

\begin{Example} \label{EX}Let $S$ be a $3$-dimensional Cohen-Macaulay local ring and $x, y, z$ \,a system of parameters
for $S$. We set $R=S/(x^2-yz)S$ and $I=(x,y)R$. Then $I$ has minimal
$j$-multiplicity. In particular, the associated graded ring ${\rm
gr}_I(R)$ is Cohen-Macaulay, and if $S$ is Gorenstein then ${\rm
gr}_I(R)$ is Gorenstein as well.
\end{Example}
\demo  Observe that $R$ is a Cohen-Macaulay local ring of dimension $d=2$.  By  \cite[4.2]{NU}, $I$ is a Cohen-Macaulay ideal of height $1=d-1$ which is generically a complete intersection and $\ell(I)=2$. Let $\xi$ be a general element of $(x,y)S$ so that $I=\xi R+yR$. By \cite[4.2]{NU},
\begin{equation}\label{un}
j(I)=\lambda_S (S/(x,y,z)S).
\end{equation}
Let $\overline{R}=R/(\xi R:_R (x,y)R^{\infty})=R/(\xi R :_R y)$, where the last equality holds since $(x,y)R_{\p}^{\infty}=yR_{\p}$ for every associated prime ideal $\p \in {\rm Ass}(R/ \xi R)$. We claim that
\begin{equation}\label{un2}
\lambda (I\overline{R}/I^2\overline{R}) = \lambda_S (S/(x,y,z)S)
\end{equation}
which together with (\ref{un}) implies that $I$ has minimal $j$-multiplicity.
To prove (\ref{un2})  we may assume that $\xi=x+\mu y$ for some $\mu \in S$. Notice that
\begin{eqnarray*}\label{un3}
\lambda (I\overline{R}/I^2\overline{R}) &=& \lambda_S( [(x,y)S+((x^2-yz, x+ \mu y)S :_S y)]/[(x^2,xy,y^2)S + ((x^2-yz, x+ \mu y )S:_S y)])\\
&=&  \lambda_S ( [(x,y)S+(x+\mu y, z + \mu x)S ]/[(x,y)^2S + (x+\mu y, z + \mu x)S ])\\
&=&  \lambda_S ( (x,y,z)S/[(x+\mu y, z + \mu x)S + (x,y,z)^2S])\\
&=&\lambda_S (S/(x,y,z)S),
\end{eqnarray*}
where the second equality follows from $(x+\mu y, z + \mu x)S =  (x^2-yz, x+ \mu y)S :_S y$. Indeed it is easy to see $(x+\mu y, z + \mu x)S \subseteq  (x^2-yz, x+ \mu y)S :_S y$. For the other inclusion, let  $\rho \in (x^2-yz, x+ \mu y) :_S y $ thus $\rho y= \eta (x+ \mu y) + \gamma(x^2-yz)$ for some elements $\eta$ and $\gamma$ of $S$, which gives $(\rho -\mu \eta + \gamma z)y=(\eta + \gamma x)x$. As $x, y$ form an $S$-regular sequence, it follows that $\rho -\mu \eta + \gamma z=\alpha x$ and $\eta + \gamma x=\alpha y$. In particular, $\rho= \alpha(x + \mu y) - \gamma(z + \mu x) \in  (x+\mu y, z + \mu x)$.
%We first show the second equality in (\ref{un3}).  which in turns yields $ (x+\mu y, z + \mu x) + (x,y,z)^2 \subseteq (x^2,xy,y^2) + ((x^2-yz, x+ \mu y) :_S y)$. To prove equality, we
%only need to show the inclusion $ ((x^2-yz, x+ \mu y) :_S y) \subseteq (x+\mu y, z + \mu x)$ as the ideal $(x^2,xy,y^2)$ is contained in both sides.
%The third equality in (\ref{un3}) holds by the proof of \cite[4.2]{NU}.
The last equality
 %(\ref{un3})
 is verified by noticing that the $S$-module  $(x,y,z)/[(x+\mu y, z + \mu x) + (x,y,z)^2]$ is cyclic with annihilator $(x,y,z)$.

 By Theorem \ref{idealminimal} and Theorem \ref{Gorenstein}, the associated graded ring ${\rm gr}_I(R)$ is Cohen-Macaulay and Gorenstein, respectively.
\QED

\ms

\begin{Example} {\rm Let $S=k[x, y, z]_{(x,y,z)}$ be a $3$-dimensional regular local ring.  Set
$R=S/(x^4-y^2z^2)S$ and $I=(x^2,y^2)R.$
%
%Let $\overline{R}=R/(x^2R:I)=R/(x^2,z^2)R=S/(x^2, z^2)S$.
%$$
%j(I)=e(I\overline{R})=\lambda(\overline{R}/y^2\overline{R})=\lambda (S/(x^2,y^2,z^2)S)=8.
%$$
%$$
%=\lambda(I\overline{R}/I^2\overline{R})=\lambda((x^2, y^2, z^2)/(x^2, y^4, z^2))=8.
%$$
Computing $j$-multiplicity as in the above example one obtains $j(I)= 8=\lambda(I\overline{R}/I^2\overline{R})$ where $\overline{R}=R/(\xi R:I)$ and $\xi$ is a general element in $I$. Now again by Example \ref{EX}  the associated graded ring ${\rm gr}_I(R)$ is Gorenstein. Indeed, ${\rm gr}_I(R)\cong {\rm gr}_{(x^2,y^2)S}(S)/((x^4-y^2z^2)^*).$}
\end{Example}

\end{document}